\documentclass[a4paper, 11pt]{article}

\usepackage[T1]{fontenc}
\usepackage[margin=1in]{geometry}
\usepackage{setspace}
\usepackage{microtype}
\usepackage[htt]{hyphenat}

\usepackage{amsmath} \allowdisplaybreaks
\usepackage{amssymb}
\usepackage{amsthm} \theoremstyle{definition}

\newtheorem{observation}{Observation}

\usepackage{natbib}
\usepackage{bibentry}

\usepackage[figuresright]{rotating}
\usepackage{algorithmic}
\usepackage{algorithm}

\usepackage{booktabs}
\usepackage{multirow}
\usepackage{graphicx}
\usepackage{appendix}
\usepackage{tabularx}

\usepackage{soul} %

\usepackage[hyphens]{url}
\usepackage[bookmarks=false,hidelinks]{hyperref}
\usepackage[capitalise,noabbrev]{cleveref}
\crefname{observation}{Observation}{Observations}
\Crefname{observation}{Observation}{Observations}
\crefname{figure}{Fig.}{Figs.}

\usepackage{tikz}
\usepackage{pgfplots}
\usetikzlibrary{positioning}
\usetikzlibrary{arrows,shapes}
\usetikzlibrary{snakes}

\usepackage{authblk}

\usepackage[shortlabels]{enumitem}

\usepackage[makeroom]{cancel}

\newcommand{\email}[1]{{\href{mailto:#1}{\nolinkurl{#1}}}}

\usepackage{bm}

\usepackage{xcolor}

\usepackage{etoolbox}
\usepackage[normalem]{ulem}
\usepackage{titlesec}
\titleclass{\subsubsubsection}{straight}[\subsubsection]
\newcounter{subsubsubsection}[subsubsection]
\renewcommand\thesubsubsubsection{\thesubsubsection.\arabic{subsubsubsection}}
\titleformat{\subsubsubsection}[hang]{\normalfont\normalsize\bfseries}{\thesubsubsubsection}{1em}{}
\titlespacing*{\subsubsubsection}{0pt}{3.25ex plus 1ex minus .2ex}{0.5em}
\makeatletter
\def\toclevel@subsubsubsection{4}
\makeatother

\makeatletter
\def\do#1{\@namedef{#1c}{\ensuremath{\mathcal{#1}}}}
\docsvlist{A,B,C,D,E,F,G,H,I,J,K,L,M,N,O,P,Q,R,S,T,U,V,W,X,Y,Z}
\makeatother

\newcommand{\dt}{\mathop{}\!\mathrm{d}t}

\renewcommand{\bar}[1]{\mkern 1.5mu\overline{\mkern-1.5mu\vphantom{x}#1\mkern-1.5mu}\mkern 1.5mu}

\usepackage{colortbl}
\usepackage{longtable}
\usepackage{pdflscape}
\usepackage{subcaption}
\usepackage{float}
\usepackage{siunitx}

\AtBeginEnvironment{algorithm}{\setstretch{1}}
\AtBeginEnvironment{table}{\setstretch{1}}
\AtBeginEnvironment{figure}{\setstretch{1}}
\AtBeginEnvironment{longtable}{\setstretch{1}}

\title{Towards LLM-Driven Automation of Lagrangian Relaxation Design for Mixed-Integer Linear Programs\thanks{A preliminary version of this work was presented at AI4Sci Korea 2026, Seoul, Korea, September 2026, and appeared in its non-archival e-proceedings \citep{sim2026ai4sci}.}}

\author[a]{Yoonju Sim\thanks{\email{syj5268@kaist.ac.kr}}}
\author[a]{Jimin Nam\thanks{\email{jimia622@kaist.ac.kr}}}
\author[a]{Jungwon Ahn\thanks{\email{ingarden117@kaist.ac.kr}}}
\author[a,b,c]{Changhyun Kwon\thanks{Corresponding author: \email{chkwon@kaist.ac.kr}}}
\affil[a]{Department of Industrial and Systems Engineering, KAIST, Daejeon 34141, Republic of Korea}
\affil[b]{Graduate School of AI for Mathematics, KAIST, Daejeon 34141, Republic of Korea}
\affil[c]{Omelet, Inc., Daejeon 34051, Republic of Korea}

\date{\today}

\begin{document}
\maketitle
\begin{abstract}
Lagrangian relaxation (LR) exploits the structure of mixed-integer linear programs (MILPs), but choosing constraints to dualize and designing primal heuristics require problem-specific expertise. Existing automated methods rely on matrix structure and problem-independent rules rather than the algebraic formulation. Large language models (LLMs), by contrast, can reason over both formulations and solver feedback. We present AutoLR, a computational laboratory where LLM agents score relaxations via Dual AUC and generate solver-free primal recovery and improvement programs. Across twenty MILP benchmarks, ten from published LR studies and ten held-out, initial LLM proposals surpass random, expert, and structure-detection baselines within three evaluations and capture most relaxation gains. Generated heuristics restore feasibility on every call across all viable benchmarks, whereas generic heuristics often fail. Lagrangian dual information offers selective, problem-dependent improvements. Together, AutoLR and its benchmark suite provide a controlled testbed for integrating LLMs into Lagrangian relaxation.
\end{abstract}


\section{Introduction}
\label{sec:introduction}

General-purpose solvers often struggle to find strong bounds or feasible solutions for mixed-integer linear programming (MILP) problems within practical time limits.
Overcoming this limitation typically requires exploiting problem structure through decomposition methods \citep{wolsey1998integer,conforti2014integer}.
Among these, Lagrangian relaxation (LR) is one of the most widely used frameworks and is closely related to Dantzig--Wolfe reformulation and column generation \citep{geoffrion1974lagrangean,fisher1981lr,frangioni2005lagrangian}.

Dualizing selected constraints into the objective with Lagrangian multipliers yields simpler subproblems whose optimal values provide valid dual bounds, while recovering feasible solutions remains a separate task.
An LR algorithm therefore relies on two key problem-specific elements, \emph{relaxation design} and \emph{primal heuristics}, which in practice are typically tailored manually to each formulation \citep{Tang2011OrderBatching,NaderiKilicDasci2020,Gao2024Lagrange}.
This manual design demands domain expertise, mathematical reasoning, and extensive computational trial-and-error, an engineering cost that limits the practical adoption of decomposition on new models.

Relaxation design faces a trade-off between bound strength and computational tractability.
Relaxing additional constraints can simplify the retained subproblems, but it may weaken the optimal dual bound and increase the number of Lagrangian multipliers \citep{geoffrion1974lagrangean,wolsey1998integer}.
Existing automated methods partition the numeric constraint matrix into independent blocks \citep{bergner2015automatic,weiner2020automatic} or score candidate relaxations with static proxies \citep{weiner2023ranking}.
Yet partitioning solely to create many small blocks can yield weak bounds or unbounded subproblems.
Conversely, retaining a large coupled block can yield a tight bound at moderate scale but slow sharply on larger instances.
Both outcomes depend on how a relaxation behaves under a finite budget, which matrix structure and static proxies cannot fully anticipate.

Even when a relaxation yields a strong dual bound, its subproblem solutions typically violate the relaxed constraints, so primal heuristics are needed to recover feasible solutions to the original model.
Effective LR heuristics repair these violations and locally improve the objective using problem-specific logic, but each is hand-crafted for a particular formulation \citep{AlaouiSelsouli2012,DingLiuZhuLi2014,GuMemarZinder2019}.
Generic alternatives such as Relaxation Induced Neighborhood Search \citep[RINS;][]{Danna2005RINS}, Relaxation Enforced Neighborhood Search \citep[RENS;][]{Berthold2014RENS}, and Feasibility Jump \citep[FJ;][]{Luteberget2023FeasibilityJump} avoid this effort by relying on problem-independent sub-MIPs or local moves.
However, they ignore what the violated constraints represent and often fail to recover feasibility within tight time limits.
Effective primal recovery thus depends on formulation-specific logic, which has so far relied largely on manual design.

For both tasks, a human expert supplies what existing automated methods lack.
The expert designs an LR algorithm for a new MILP by reading the formulation, reasoning about the resulting subproblem structures, evaluating promising relaxations computationally, and crafting primal heuristics from observed violations.
Automating this process requires a distinct capability for each step: algebraic reasoning over the formulation to identify constraint roles and select which to relax, interpretation of numerical feedback under finite time budgets, and problem-specific code generation to recover feasible solutions.
Because these capabilities must work together in a sequential workflow, LR offers a natural testbed for \emph{agentic} operations research.

Recent large language model (LLM) research in optimization has progressed along two main fronts: automating mathematical modeling from natural language \citep{ahmaditeshnizi2024optimus,huang2025orlm} and designing heuristics or solver components \citep{novikov2025alphaevolve,baek2026llms,kong2026frontieror,cheng2026decompgym,hou2026llm4branch}.
Yet automated modeling hands the resulting MILP directly to a general-purpose solver, heuristic and end-to-end algorithm frameworks seek primal solutions without dual bounds, and existing exact-method agents either tune low-level solver rules or implement human-prescribed decompositions rather than deciding how to decompose a new model.
Automating both constraint selection and paired primal recovery directly from an algebraic MILP therefore remains largely unexplored.

To close this gap, we present \emph{AutoLR}, a computational laboratory in which two LLM agents design, test, and refine LR algorithms over a unified solver engine.
The Relaxation Design Agent proposes constraint selections from the formulation, tests them under a finite solver budget, and refines them from solver feedback before selecting one relaxation for the model.
Each candidate is scored by \emph{Dual AUC}, a time-integrated measure of the normalized dual-bound trajectory over the budget.
Given the selected relaxation, the Primal Heuristic Agent writes, tests, and refines solver-free programs for initial recovery and improvement.
The resulting relaxation and executable heuristic are fixed offline at the formulation level, so online solving requires no LLM calls.
Because AutoLR organizes this workflow into modular stages rather than a single black box, it also lets us isolate the contributions of initial proposals, iterative search, and Lagrangian dual information against controlled baselines.
To evaluate this workflow under finite solver budgets and test generalization beyond published LR studies, we pair Dual AUC with a twenty-benchmark MILP suite combining ten published LR formulations with ten held-out MILPs that have no published LR designs.

Across these twenty benchmarks, our experiments clarify both the strengths and the limits of LLM-guided decomposition.
When paired with solver evaluation, reading the algebraic formulation enables both frontier and compact open-weight models to propose initial relaxation pools that outperform random and structure-detection baselines.
Subsequent iterative search, whether LLM-guided or based on Bayesian optimization (BO), mainly improves run-to-run consistency.
In primal heuristic generation, formulation-aware code synthesis recovers feasible solutions on unseen hard instances for every benchmark whose selected relaxation remains viable at hard scale, whereas generic MILP heuristics frequently fail.
Lagrangian dual information yields selective, problem-dependent gains in primal recovery.
Most designs obtained on medium-scale instances also transfer to harder unseen instances, although subproblem scaling occasionally becomes a bottleneck.
Together, these results suggest that offline, formulation-level automation can reduce the engineering cost of applying decomposition to new models.

This paper makes three main contributions.
First, we devise AutoLR, an agentic workflow that designs reusable constraint-family Lagrangian relaxations and solver-free primal heuristics at the formulation level, so that all LLM costs are incurred offline and online solving requires no LLM calls.
Second, we propose the Dual AUC metric for evaluating finite-budget dual-bound progress and release a twenty-benchmark MILP suite spanning seven application domains, which separates published LR formulations from held-out MILPs without published LR designs.
Third, we evaluate relaxation design across frontier and compact open-weight models against classical search and structure-detection baselines, and evaluate primal recovery with and without Lagrangian dual information against generic MILP heuristics.
Together, these experiments separate the roles of formulation-based initialization and iterative search, and characterize when medium-scale designs transfer to harder instances.
All twenty benchmark formulations, baseline relaxations, instance-level analyses, code, and data are available at \url{https://github.com/syj5268/AutoLR}.

The remainder of this paper is organized as follows.
\Cref{sec:background} provides background on LR and reviews related work, and \cref{sec:methodology} presents AutoLR.
\Cref{sec:experimental_results} reports the main results for relaxation design and primal heuristic generation.
\Cref{sec:instance_wise_results} examines selected relaxation structures, scale transfer, search trajectories, and joint primal--dual progress.
Finally, \cref{sec:discussion} discusses implications and limitations and concludes the paper.

\section{Background and Related Work}
\label{sec:background}
This section introduces LR and its notation, then reviews related work on constraint selection, primal heuristics, and LLMs in optimization.

\subsection{Lagrangian Relaxation}
\label{sec:background_lr}
Classical decomposition approaches for MILPs include Benders decomposition, Dantzig--Wolfe reformulation with column generation, and Lagrangian relaxation \citep{wolsey1998integer}.
Among them, LR moves complicating constraints directly into the objective with Lagrangian multipliers to obtain valid dual bounds and simpler subproblems \citep{geoffrion1974lagrangean,fisher1981lr}.
Consider an MILP $\min_{x\in X}\{c^\top x\mid Ax\geq b\}$ with finite optimal value $v^\ast$, where $X$ imposes variable bounds and integrality requirements. Equality constraints follow the same formulation with unrestricted multipliers.
Partitioning the rows into relaxed and retained sets $(A_R,b_R)$ and $(A_F,b_F)$ defines $X_F=\{x\in X\mid A_Fx\geq b_F\}$.
In practice, this partition operates on \emph{constraint families}, indexed sets of constraints sharing a common algebraic expression, so each family is relaxed or retained as a whole.
For multipliers $\lambda\geq0$, the Lagrangian dual function is
\begin{equation}
L(\lambda)=\min_{x\in X_F}\bigl\{c^\top x+\lambda^\top(b_R-A_Rx)\bigr\},
\label{eq:lagrangian}
\end{equation}
and the Lagrangian dual is $v_{\mathrm{LD}}=\sup_{\lambda\geq0}L(\lambda)\leq v^\ast$.
When $X_F$ separates into variable blocks, the relaxed problem \eqref{eq:lagrangian} decomposes into independent subproblems, though even a single retained block can yield a useful bound if it is easier to solve than the original model.

The tightness of this bound is determined by the convex hull of $X_F$ \citep{geoffrion1974lagrangean}, since
\begin{equation}
v_{\mathrm{LD}}=\min\{c^\top x\mid A_Rx\geq b_R,\ x\in\operatorname{conv}(X_F)\}.
\label{eq:lagrangian_convexification}
\end{equation}
If $X_F$ has the integrality property, $\operatorname{conv}(X_F)$ equals its continuous relaxation and $v_{\mathrm{LD}}$ reduces to the LP relaxation value, whereas tighter bounds arise when $\operatorname{conv}(X_F)$ strengthens the LP relaxation while keeping $X_F$ tractable.
For the same retained set $X_F$, $v_{\mathrm{LD}}$ also equals the LP bound of the corresponding Dantzig--Wolfe reformulation, which solves the same subproblems during column generation.

\begin{figure}
  \centering
  \includegraphics[width=\textwidth]{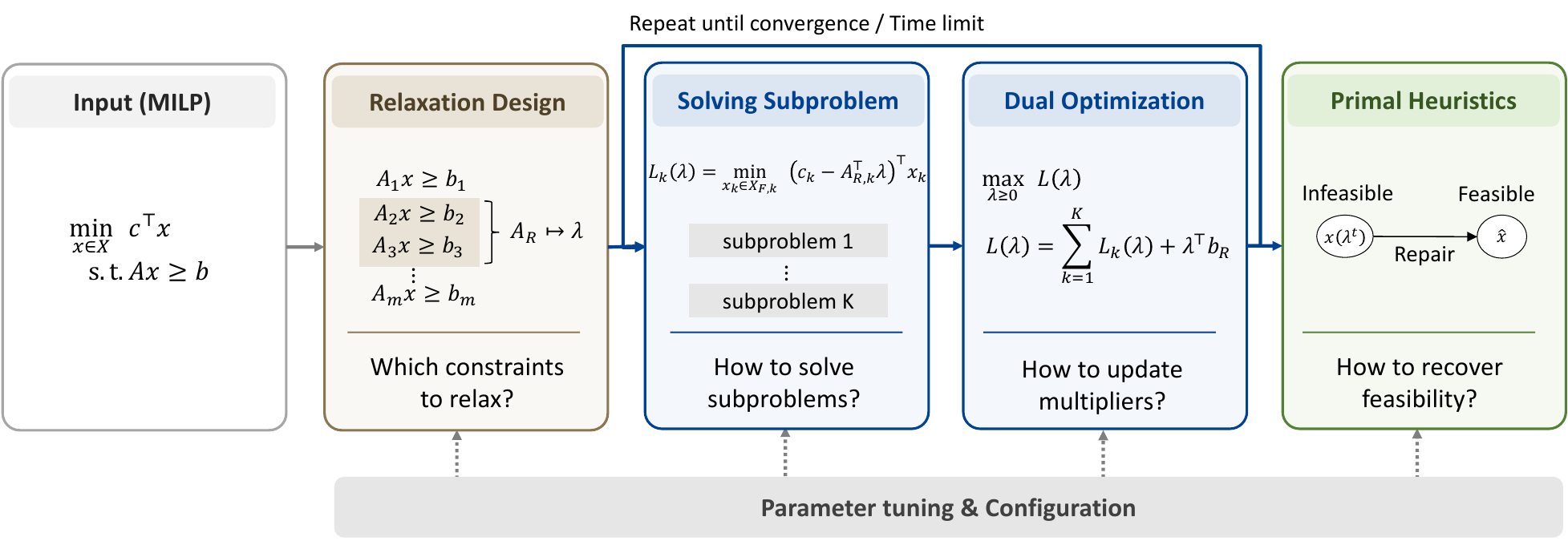}
  \caption{Workflow of a Lagrangian relaxation algorithm and its core design components.}
  \label{fig:design_space}
\end{figure}

As summarized in \cref{fig:design_space}, building an LR algorithm around a chosen relaxation requires a subproblem solver, a dual optimization method, and a primal heuristic.
Subproblems can be solved by specialized algorithms for known structures or by a generic MILP solver.
Because $L(\lambda)$ is concave and nonsmooth, dual optimization relies on subgradient, volume, or proximal bundle methods \citep{fisher1981lr,kiwiel1990proximity,feltenmark2000dual}.
Among them, bundle and volume methods aggregate subproblem solutions $x(\lambda)$ across iterations into a point $\bar x_t\in\operatorname{conv}(X_F)$ \citep{barahona2000volume,frangioni2005lagrangian}.
Because subproblem solutions $x(\lambda)$ generally violate $A_Rx\geq b_R$, and even $\bar x_t$ may be fractional or infeasible, primal heuristics are needed to recover feasible solutions to the original model.
Whereas generic MILP solvers and bundle methods can be reused across formulations, relaxation design and primal heuristics remain highly problem-dependent.

\subsection{Relaxation Design and Primal Heuristics}
\label{sec:background_relaxation}\label{sec:background_recovery}
In classical LR applications, both the choice of relaxed constraints and the primal heuristic are crafted manually for each problem class \citep{fisher1981lr,frangioni2005lagrangian}.
This subsection reviews classical practice and existing automation attempts for relaxation design on the dual side and primal heuristics on the primal side.

On the dual side, selecting which constraint families to relax affects bound tightness, subproblem tractability, and the difficulty of dual optimization.
Published LR studies typically make this choice to preserve tractable structure in the retained problem.
For example, \citet{NaderiKilicDasci2020} retain knapsack and assignment subproblems in a capacitated transshipment problem, and \citet{Gao2024Lagrange} relax capacity constraints to obtain single-item lot-sizing subproblems.
Without domain knowledge, automatic structure-detection methods identify linking rows that separate the constraint matrix into independent blocks and score partitions using structural metrics such as relative border area and block balance \citep{bergner2015automatic,Khaniyev2018StructureDetection} or data-driven models \citep{basso2023datadriven}.
However, these methods favor block-diagonal partitions without considering dual-bound quality or convergence.
\begin{figure}
  \centering
  \begin{subfigure}[t]{0.42\textwidth}
    \centering
    \includegraphics[width=\linewidth]{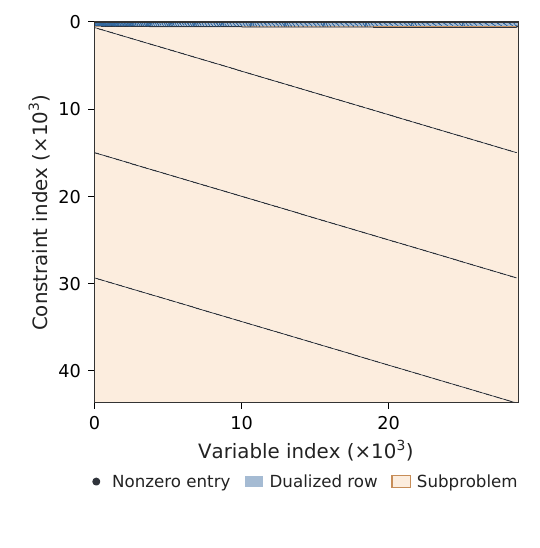}
    \caption{Expert: 500 relaxed rows, 1 block; bound gap = 5.3\%.}
    \label{fig:structure_example_expert}
  \end{subfigure}\hspace{0.04\textwidth}
  \begin{subfigure}[t]{0.42\textwidth}
    \centering
    \includegraphics[width=\linewidth]{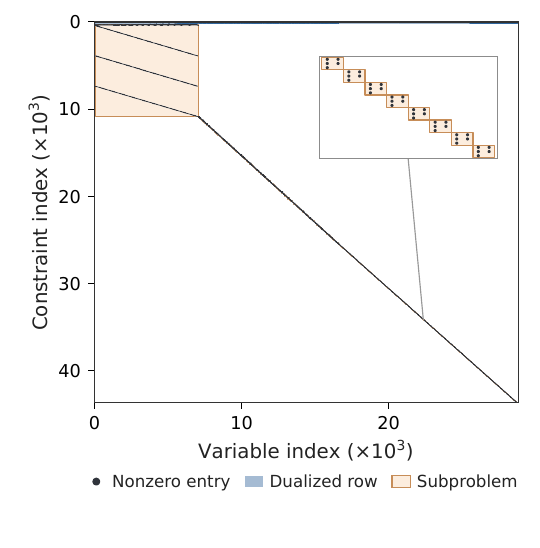}
    \caption{hMETIS: 205 relaxed rows, 6,356 blocks; bound gap $> 100\%$.}
    \label{fig:structure_example_hmetis}
  \end{subfigure}
  \caption{Constraint--variable incidence structures of two relaxations on benchmark P05 \citep{Tang2011OrderBatching}. The inset in (b) magnifies the diagonal subproblem blocks. Bound gap is $100(U-b)/|U|$ within 120\,s, where $b$ is the best LR lower bound and $U$ is the Gurobi incumbent.}
  \label{fig:structure_example}
\end{figure}
For example, in \cref{fig:structure_example}, the single-block relaxation of \citet{Tang2011OrderBatching} for a steelmaking order-batching problem achieves a small bound gap, whereas hMETIS partitioning \citep{bergner2015automatic}, which splits the rows into balanced blocks while minimizing the number of linking rows, creates thousands of blocks yet encounters an unbounded subproblem during dual optimization, leaving a negative best bound and a bound gap exceeding 100\%.
Even among methods designed for LR, \citet{weiner2020automatic} search over relaxations using matrix structure alone, and \citet{weiner2023ranking} rank candidates using structural features and a single dual-function evaluation, not an iterative dual-bound trajectory.

On the primal side, an LR heuristic repairs violations of the relaxed constraints and locally improves a feasible solution $\widehat x\in X$ satisfying $A\widehat x\geq b$ using subproblem solutions, Lagrangian multipliers, and aggregate primal points \citep{fisher1981lr,frangioni2005lagrangian}.
Without problem-specific design, generic MILP heuristics rely on problem-independent strategies such as the Feasibility Pump \citep{fischetti2005feasibility} and neighborhood sub-MIP heuristics including RINS \citep{Danna2005RINS} and RENS \citep{Berthold2014RENS}.
Aggregate primal points can guide such neighborhood sub-MIPs, although success depends on the fractional point and the time limit on the sub-MIP.
Alternatively, FJ avoids solving sub-MIPs and uses local search over variable values to reduce weighted constraint violations \citep{Luteberget2023FeasibilityJump}.
Generic heuristics thus sacrifice formulation-specific repair logic, whereas problem-specific Lagrangian heuristics require manual design.

\subsection{LLMs in Optimization}
\label{sec:related_work}
LLM research in optimization ranges from automated formulation \citep{ahmaditeshnizi2024optimus,huang2025orlm} to automated algorithm design, spanning both heuristic and exact methods.
For heuristic algorithm design, evolutionary frameworks such as AlphaEvolve and ReEvo evolve heuristic code \citep{novikov2025alphaevolve,ye2024reevo}, and AutoCO coevolves temporary constraint-relaxation heuristics for primal search \citep{liu2026autoco}.
QuantumEvo similarly evolves reusable variable-ordering heuristics for binary decision diagrams, evaluating them by the quantum cost of synthesized circuits \citep{sim2026quantumevo}.
Like AutoLR's primal heuristics, the discovered programs run without online LLM calls.
With advances in frontier LLMs, greedy refinement \citep{sun2025cobench} and even single-query generation \citep{baek2026llms} also produce strong primal heuristics.
For large-scale OR problems drawn from the literature, FrontierOR benchmarks end-to-end algorithm code generation by feasible solution quality and runtime against Gurobi, finding that designing scalable algorithms beyond direct solver calls remains difficult \citep{kong2026frontieror}.
For exact algorithm design, LLMs generate branching and cutting-plane strategies \citep{hou2026llm4branch,zheng2026cutting}, while DecompGym tests coding agents that implement predefined Benders or column-generation decompositions \citep{cheng2026decompgym}.

When candidate evaluations are costly, LLM-guided algorithm design also intersects with black-box optimization.
Although LLMs can warm-start BO or guide sequential search \citep{liu2024llambo,roohani2025biodiscovery}, recent benchmarks find that iterative LLM search often matches rather than outperforms classical BO \citep{gupta2025llmbo,demir2026finitepool}.

Building on this literature, \cref{sec:methodology} introduces AutoLR.
Unlike systems that implement a predefined decomposition or search only for primal incumbents, AutoLR selects constraint families, separates initialization from iterative search, and generates primal heuristics tailored to the chosen relaxation.

\section{AutoLR: A Computational Laboratory for Lagrangian Relaxation}
\label{sec:methodology}

\begin{figure}
  \centering
  \includegraphics[width=\textwidth]{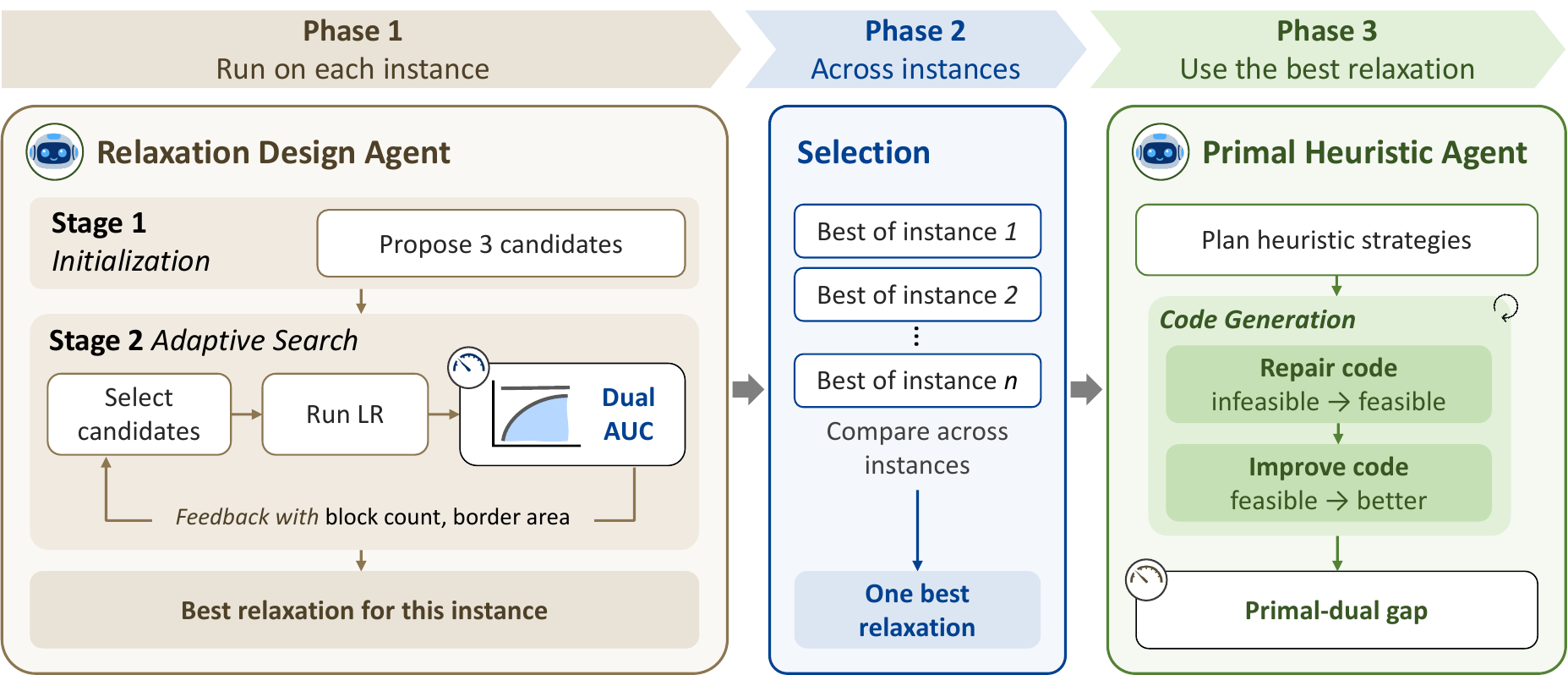}
  \caption{AutoLR framework: a unified computational environment with two agents for relaxation design and primal heuristic generation.}
  \label{fig:framework}
\end{figure}

We present AutoLR, a computational laboratory for designing and evaluating Lagrangian relaxation algorithms, as illustrated in \cref{fig:framework}.
It combines a unified computational environment with two LLM-driven agents, the Relaxation Design Agent and the Primal Heuristic Agent.
The environment provides a modular LR engine, structural analysis, and evaluation procedures while holding the subproblem solver and dual optimization method fixed across all comparisons.

Given a formulation $F$ and a set of instances $\mathcal I$, AutoLR integrates both agents into a three-phase agentic workflow.
Phase~1 searches over constraint-family relaxations on each instance in $\mathcal I$, Phase~2 selects a single relaxation $z^\ast$ for the formulation, and Phase~3 generates a repair-and-improvement program $P^\ast$ for $z^\ast$.
Together with the LR engine, $(z^\ast,P^\ast)$ forms an LR algorithm that can be applied to unseen instances without further design. 
\Cref{sec:relaxation_design_agent} details the evaluation metric with Phases~1 and~2, \cref{sec:primal_recovery_agent} details Phase~3, and Algorithm~\ref{alg:autolr} in Appendix~\ref{app:autolr_alg} summarizes the complete three-phase workflow.

\subsection{Relaxation Design Agent}
\label{sec:relaxation_design_agent}
The Relaxation Design Agent seeks a relaxation with high dual-bound performance on instance $I$ within an LR time budget $T$.
It selects the constraint families to dualize.
Let $m_{\mathrm{model}}$ denote the number of constraint families in the formulation.
A candidate relaxation is represented by a binary mask $z\in\{0,1\}^{m_{\mathrm{model}}}$, where $z_j=1$ dualizes every row in constraint family $j$ and $z_j=0$ retains it.
The remainder of this subsection describes candidate evaluation, iterative search, and final relaxation selection.

\subsubsection{Evaluation by Dual AUC}
\label{sec:autolr_evaluator}
Evaluating candidate relaxations requires a metric of their dual-bound performance under a fixed computational budget.
Let $v_{\mathrm{LD}}(z)$ denote the optimal Lagrangian dual value for mask $z$.
For a bound $b$ obtained within $T$, the gap $v^\ast-b$ splits into the Lagrangian duality gap $v^\ast-v_{\mathrm{LD}}(z)$ and the remaining gap $v_{\mathrm{LD}}(z)-b$.
The duality gap is determined by the dualized families alone, whereas the remaining gap also depends on the computation completed within $T$.
As an immediate consequence of \eqref{eq:lagrangian_convexification} \citep{geoffrion1974lagrangean}, the following observation characterizes how additional dualization affects $v_{\mathrm{LD}}(z)$.

\begin{observation}
\label{obs:monotone}
For minimization, $z'\le z$ implies $v_{\mathrm{LD}}(z)\le v_{\mathrm{LD}}(z')\le v^\ast$.
\end{observation}

\Cref{obs:monotone} shows that keeping more constraints can strengthen the optimal dual bound.
These constraints may also make the subproblems harder to solve and thus slow the approach to that bound.
Both effects determine the bound reached within $T$, but the final bound alone depends on the particular choice of $T$ and does not reflect when that bound was reached.

Therefore, building on \citeauthor{berthold2013measuring}'s (\citeyear{berthold2013measuring}) discussion of primal and dual integrals, we define \emph{Dual AUC}, denoted by $S(z;I,T)$, as the average normalized dual-bound quality over $[0,T]$:
\begin{equation}
S(z;I,T)=\frac{1}{T}\int_0^T q_I\bigl(b_z(t)\bigr)\, \dt
=\mathbb{E}_{\tau\sim\mathcal{U}(0,T)}\!\left[q_I\bigl(b_z(\tau)\bigr)\right],
\label{eq:dual_auc}
\end{equation}
where $b_z(t)$ is the best dual bound available for candidate $z$ on instance $I$ at time $t$, and $q_I(b)$ is its normalized quality against a fixed feasible reference value $U_I>0$.
For minimization, $q_I(b)=\max\{0,b/U_I\}$, whereas for maximization, $q_I(b)=U_I/b$ when $b>0$ and zero otherwise.
If $U_I<0$, negating the objective converts between minimization and maximization so that $U_I>0$.
The integrand is zero before the first finite bound, and bounds obtained after $T$ contribute no area.
Equivalently, Dual AUC is the expected normalized quality of the best available bound at a time drawn uniformly from $[0,T]$.
This interpretation gives equal weight to all times in $[0,T]$ without assuming that computational budgets in practice are uniformly distributed.

\subsubsection{Relaxation Search and Selection}
\label{sec:design_and_search}
To obtain a single relaxation $z^\ast$ for formulation $F$, Phase~1 searches over candidate masks on each instance $I\in\mathcal I$ across $R$ replicates, and Phase~2 consolidates those runs into one selection.
Phase~1 splits each search run into two stages, initial candidate generation and iterative refinement, so that formulation-based reasoning can be distinguished from feedback-driven search.
The Relaxation Design Agent uses LLMs in both stages. 
The same two-stage protocol also supports random initialization and random or combinatorial BO search~\citep{baptista2018bocs} as controlled baselines.

\label{sec:initial_design}
In Phase~1, Stage~1 generates the initial candidate pool before any LR run to isolate what the agent infers from the objective, indexed constraint expressions, and family row counts.
For each instance $I\in\mathcal I$ and replicate $r=1,\ldots, R$, the agent proposes three distinct masks that are evaluated under budget $T$ to seed the search history $\mathcal H$.
The prompt asks for relaxations that reduce subproblem cost while preserving bound strength.

Then, Stage~2 refines the pool iteratively.
In each round, the agent first proposes a shortlist of unevaluated masks.
Before running LR, a structural preview $\pi(c;I)$ is computed for each shortlisted candidate $c$ without solving subproblems, drawing on the same structural descriptors used in automatic structure detection \citep{bergner2015automatic}: the connected-component block count and the relative border area, namely the fraction of dualized rows.
The selection step then inspects these previews alongside $\mathcal H$ to choose one mask to evaluate by Dual AUC under budget $T$ and append to $\mathcal H$.
Proposal rationales are omitted so that selection relies solely on structural previews and observed Dual AUC, avoiding verbosity bias \citep{zheng2023judging}.
Search continues until an evaluation cap $N$ or stopping rule is met and returns the highest-scoring mask $z_{I,r}$.

\label{sec:autolr_selection}
Phase~2 consolidates repeated searches across $\mathcal I$ into a single relaxation $z^\ast$ to limit sensitivity to instance-specific and single-run variation.
If no run achieves a positive Dual AUC on $\mathcal I$, Phase~2 selects no relaxation.
When a single mask wins strictly more than half of the runs, Phase~2 selects that majority winner.
Otherwise, it evaluates any missing $(z,I)$ pairs among the distinct run winners and selects the mask with the highest mean Dual AUC over $\mathcal I$.

\subsection{Primal Heuristic Agent}
\label{sec:primal_recovery_agent}
The Primal Heuristic Agent constructs a solver-free program $P=(R_\phi,Q_\psi)$ in Phase~3 for the selected relaxation $z^\ast$.
This subsection details the program structure, gap-based evaluation and verification, and the plan-and-revise loop.
The program separates initial recovery $R_\phi$ from improvement $Q_\psi$ so that a feasible incumbent is established before anytime local search begins.
Each candidate receives the instance's variables, objective coefficients, constraints, and a frozen terminal state $s_I=(x^{\mathrm{LR}},\lambda,\bar x,\mathcal X)$ containing the final subproblem solution, Lagrangian multipliers, bundle aggregate primal point, and past iteration records of $(x^{\mathrm{LR}},\lambda,L(\lambda))$.
The initial recovery function $R_\phi(I,z^\ast,s_I)$ returns a first feasible solution $\widehat x$ from $s_I$ without restricting its internal strategy.
When $\widehat x$ is feasible, the improvement function computes $\widetilde x=Q_\psi(I,z^\ast,s_I,\widehat x)$ and replaces $\widehat x$ only if $\widetilde x$ is feasible and strictly better.
Both functions share a single time limit per call, and generated programs cannot call external solvers.

Each generated program $P$ is evaluated on an instance $I$ by its primal--dual gap against the fixed LR bound $d_I$,
\begin{equation}
\gamma_I(P)=\frac{|f_I(P)-d_I|}{\max\{|f_I(P)|,\epsilon\}},
\label{eq:mean_lr_gap}
\end{equation}
where $f_I(P)$ is the recomputed objective of a feasible solution, $\epsilon>0$ guards against zero denominators, and any call that fails to return a feasible solution within the time limit is logged with $\gamma_I(P)=1$ and rejects the candidate. A verifier checks that $P$ neither calls external solvers nor modifies the model, runs it in an isolated subprocess, validates all variable domains and constraints, and recomputes $f_I(P)$.
Candidate programs are then ranked by the mean gap $\bar\gamma(P)=\frac{1}{|\mathcal I|}\sum_{I\in\mathcal I}\gamma_I(P)$ over $\mathcal I$.

The verifier also checks that the initial recovery function uses the dual state $s_I$.
It replays $R_\phi$ with the dual state removed and rejects any program whose output remains identical across every instance in $\mathcal I$.
Because using $s_I$ does not guarantee better objectives, \cref{sec:primal_recovery_results} compares this agent against an LR-blind variant that constructs an initial feasible solution directly from the instance without $s_I$.

The agent uses a plan-and-revise loop rather than population-based evolution, motivated by recent studies on few-query synthesis and greedy refinement \citep{sun2025cobench,baek2026llms}.
It first plans and implements three distinct strategies, allowing at most five error-driven corrections per strategy to fix syntax or runtime errors.
A candidate passes only if it returns a feasible solution on every instance in $\mathcal I$, and the lowest-gap program becomes $P^\ast$.
Three sequential revision rounds then use per-instance objectives, gaps, runtimes, and failure logs to update $P^\ast$ whenever a revision lowers $\bar\gamma$.
The final program $P^\ast$ is paired with $z^\ast$ and can be applied to unseen instances.

\section{Experimental Results}
\label{sec:experimental_results}\label{sec:relaxation_design_results}
This section evaluates AutoLR across twenty MILP benchmarks.
After summarizing the experimental setup in \cref{sec:experimental_setup}, the experiments address three questions: how much of the relaxation gain comes from formulation-based initialization and iterative search in \cref{sec:initialization}, how LLM choice and pool size affect initialization in \cref{sec:llm_choice}, and how generated solver-free heuristics perform on unseen hard instances and what Lagrangian dual states contribute in \cref{sec:primal_recovery_results}.

\subsection{Setup}
\label{sec:experimental_setup}\label{sec:benchmarks}\label{sec:policies}\label{sec:pr_setup}
All experiments run on an AMD Ryzen~9 5900X workstation with 24 CPUs using Julia, JuMP~1.29, and Gurobi~13.0.
The LR engine uses Gurobi as the subproblem solver and a proximal bundle method initialized at $\lambda=0$ with a $10^{-5}$ optimality tolerance under a budget of $T=120$ seconds, and primal solutions are verified to a $10^{-6}$ feasibility tolerance.
Unless noted otherwise, both agents use Claude Opus~5.5 with high reasoning effort, and token usage and API costs are detailed in Appendix~\ref{app:short_llm_cost}.
Agents have no internet access and no tools such as web search, browsing, or code execution. Each model is queried as a single structured-output API request: Claude Opus~5.5 via the Anthropic Messages API, GPT-5.6-sol via the OpenAI Chat Completions API, and Qwen3.8-27B on a self-hosted vLLM server.

Each benchmark comprises three medium and three hard instances, with complete formulations and data in the repository.
Medium instances are used for algorithm design, whereas hard instances, which Gurobi does not solve to optimality within 120 seconds, are used to evaluate transfer.
The validation set, P01--P10, reconstructs MILPs from published LR studies with their \emph{expert} relaxations, and the held-out set, T01--T10, consists of ten MILPs without a published LR design.
Appendices~\ref{app:short_validation} and~\ref{app:short_heldout} summarize both benchmark suites, including their sources, Gurobi MIP gaps, and relaxations.
The $T$-second Gurobi incumbent serves as the reference value $U_I$.

\subsection{Relaxation Design: Initialization and Search}
\label{sec:initialization}\label{sec:search}

\begin{table}
  \centering
  \caption{Initialization and search performance on the validation and held-out sets. Entries report benchmark-averaged Dual AUC $\pm$ standard deviation across three replicates, where Pool best is the best initial candidate and Final is the best after search. Best (\%) is the percentage of runs reaching the highest final Dual AUC observed across the search policies on each instance.}
  \label{tab:main_stage_results}\label{tab:heldout_main_results}
  \footnotesize
  \setlength{\tabcolsep}{3.5pt}
  \begin{tabular}{@{}ll ccr ccr@{}}
    \toprule
    & & \multicolumn{3}{c}{\textbf{Validation}} & \multicolumn{3}{c}{\textbf{Held-out}} \\
    \cmidrule(lr){3-5}\cmidrule(l){6-8}
    \textbf{Init.} & \textbf{Search} & \textbf{Pool best} & \textbf{Final} & \textbf{Best (\%)} & \textbf{Pool best} & \textbf{Final} & \textbf{Best (\%)} \\
    \midrule
    Expert & None & \multicolumn{2}{c}{0.675} & --- & \multicolumn{2}{c}{---} & --- \\
    hMETIS & None & \multicolumn{2}{c}{0.540} & --- & \multicolumn{2}{c}{0.379} & --- \\
    \midrule
    \multirow{2}{*}{Random} & Random & \multirow{2}{*}{0.548\,$\pm$\,0.024} & 0.791\,$\pm$\,0.011 & 18.3 & \multirow{2}{*}{0.426\,$\pm$\,0.021} & 0.503\,$\pm$\,0.020 & 17.2 \\
    & BO & & 0.743\,$\pm$\,0.014 & 37.2 & & 0.542\,$\pm$\,0.008 & 39.4 \\
    \midrule
    \multirow{2}{*}{LLM} & BO & \multirow{2}{*}{\textbf{0.808}\,$\pm$\,0.003} & 0.835\,$\pm$\,0.002 & 53.9 & \multirow{2}{*}{\textbf{0.512}\,$\pm$\,0.007} & 0.550\,$\pm$\,0.007 & 35.0 \\
    & LLM & & \textbf{0.840}\,$\pm$\,0.001 & \textbf{84.4} & & \textbf{0.556}\,$\pm$\,0.021 & \textbf{43.9} \\
    \bottomrule
  \end{tabular}
\end{table}

\begin{figure}
  \centering
  \captionsetup[subfigure]{font=small}
  \begin{subfigure}[t]{0.49\textwidth}
    \centering
    \includegraphics[width=\linewidth]{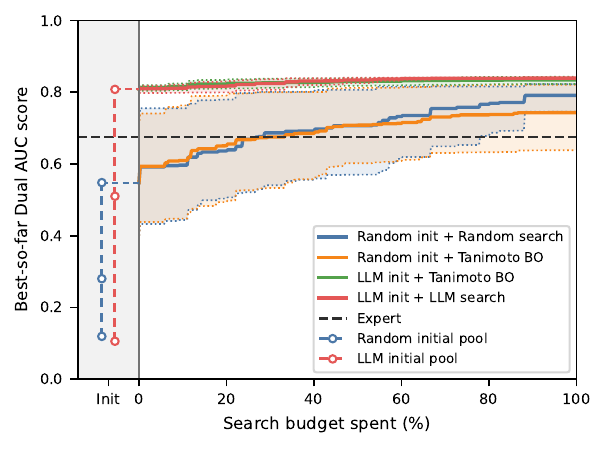}
    \caption{Validation set}
    \label{fig:aggregate_dual_auc_curve}
  \end{subfigure}\hfill
  \begin{subfigure}[t]{0.49\textwidth}
    \centering
    \includegraphics[width=\linewidth]{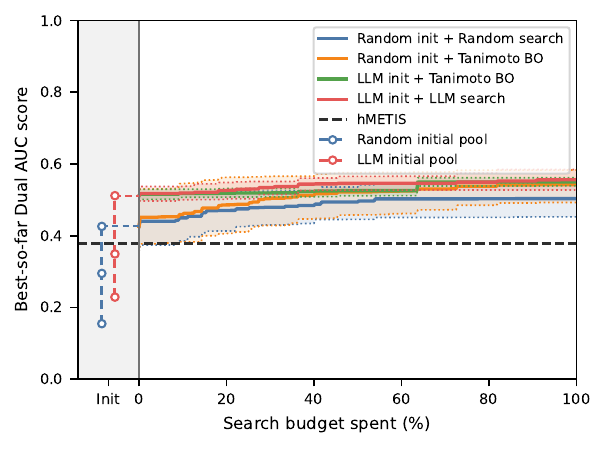}
    \caption{Held-out set}
    \label{fig:heldout_aggregate_curve}
  \end{subfigure}
  \caption{Best-so-far Dual AUC during relaxation search. Curves show benchmark-averaged scores, shaded bands span replicate minima and maxima, and Init shows the three initial scores sorted within each run before averaging across benchmarks.}
  \label{fig:aggregate_dual_auc}
\end{figure}

To separate formulation-based initialization from iterative search, Phase~1 is evaluated independently on all six instances of each benchmark.
Each run searches nonempty constraint-family masks with at most $k_{\max}=\lfloor m_{\mathrm{model}}/2\rfloor+1$ dualized families across $R=3$ replicates and up to $N=2m_{\mathrm{model}}$ evaluations per instance.
On T06--T10, which include benchmarks with larger search spaces, runs also stop early after five consecutive rounds with less than 1\% relative Dual AUC gain.
We evaluate four initialization--search policies sharing three-candidate initial pools: Random--Random, Random--BO, LLM--BO, and LLM--LLM.
They are compared by Dual AUC against the expert relaxation and the hMETIS partition over $k\in\{2,\dots,8\}$ blocks with the fewest linking rows, with BO detailed in Appendix~\ref{app:short_bo}.
Dual AUC excludes LLM calls and startup time, and its sensitivity to the time budget $T$ is analyzed in Appendix~\ref{sec:time_budget_sensitivity}.

As shown in \cref{tab:main_stage_results,fig:aggregate_dual_auc}, formulation-based initialization drives most of the relaxation gain before any computational feedback.
Reading only the algebraic objective, indexed constraints, and row counts, the three-candidate LLM pool outperforms random pools on both suites and significantly exceeds the expert relaxation at the 5\% level under two-sided Wilcoxon signed-rank tests over benchmark means.
Its mean Dual AUC also exceeds hMETIS on both suites, and after Stage~2 search LLM--LLM surpasses hMETIS significantly on both suites at the 5\% level.
On most validation benchmarks, the initial pool recovers the published expert relaxation and also proposes alternatives better suited to the generic subproblem solver.
Although validation scores may partly reflect exposure to published LR studies, the gains on the held-out set show that LLMs propose effective relaxations even without published LR designs.
Moreover, the wide spread across the three initial proposals in \cref{fig:aggregate_dual_auc} shows that even on the validation set, the LLM acts as a hypothesis generator that requires solver evaluation rather than a reliable lookup of known relaxations.

Iterative search from a strong LLM pool adds consistency rather than large mean gains, matching combinatorial BO on average.
In \cref{tab:main_stage_results}, LLM--LLM and LLM--BO raise mean Dual AUC by similar small amounts on both suites, and paired Wilcoxon tests do not separate the two search policies.
Clearer differences appear in the Best (\%) column of \cref{tab:main_stage_results}.
Whereas the initial LLM pool contains the best observed relaxation in 21.7\% of validation runs and 20.6\% of held-out runs, Stage~2 search substantially increases this discovery rate.
LLM--LLM reaches the best relaxation most consistently on both suites with the smallest replicate variation on validation, while its margin over LLM--BO narrows on held-out problems where occasional hard runs still stall at zero.

Without LLM initialization, combinatorial BO is much less sample-efficient under the short evaluation budget.
On the validation set, Random--BO finishes well below both LLM-initialized policies and even below Random--Random, with a significant paired difference from LLM--BO at $p=0.002$.
Because the Tanimoto kernel captures only mask overlap without family sizes or incidence structure, a weak or zero-scoring random pool can trap BO in local neighborhoods around poor initial masks under the short budget, whereas random search explores disjoint masks more broadly.
On the held-out set, Random--BO closes most of its initial deficit by the end of the budget, though \cref{fig:aggregate_dual_auc} shows that surrogate search spends most of its evaluations approaching where LLM pools begin.

The difference in absolute Dual AUC between the two suites partly reflects benchmark composition.
Whereas the validation benchmarks originate from published LR studies, the held-out benchmarks have no published LR designs and include formulations on which both AutoLR and hMETIS perform poorly under the current solver and time budget.
In particular, T03, T08, and T09 have mean medium-instance Dual AUC below 0.1 for every evaluated relaxation from AutoLR and hMETIS, consistent with the structural obstacles examined in \cref{sec:selected_relaxations}.
As a post hoc descriptive analysis, excluding these three benchmarks raises the final LLM--LLM mean Dual AUC across medium and hard instances from 0.556 to approximately 0.753, compared with 0.840 on validation.
This concentration of low scores suggests that the aggregate difference partly reflects the inclusion of formulations that are challenging for the evaluated LR designs.
Beyond these aggregate results, benchmark-level relaxation structures, Phase~2 selections, and medium-to-hard transfer are analyzed in \cref{sec:selected_relaxations}.

\subsection{Effect of LLM Choice and Pool Size}
\label{sec:llm_choice}

\begin{table}
  \centering
  \caption{Effect of LLM choice and initial pool size, aggregated as in \cref{tab:main_stage_results}. Dual AUC reports the best candidate in a pool of $k$ proposals, and Zero (\%) is the percentage of individual candidates scoring zero across both suites. Bold marks the best model in each column.}
  \label{tab:llm_init}
  \footnotesize
  \setlength{\tabcolsep}{4.5pt}
  \begin{tabular}{@{}l cc cc cc@{}}
    \toprule
    & \multicolumn{2}{c}{\textbf{Validation}} & \multicolumn{2}{c}{\textbf{Held-out}} & \multicolumn{2}{c}{\textbf{Zero (\%)}} \\
    \cmidrule(lr){2-3}\cmidrule(lr){4-5}\cmidrule(l){6-7}
    \textbf{Model} & $\boldsymbol{k=1}$ & $\boldsymbol{k=3}$ & $\boldsymbol{k=1}$ & $\boldsymbol{k=3}$ & $\boldsymbol{k=1}$ & $\boldsymbol{k=3}$ \\
    \midrule
    Claude Opus~5.5 & \textbf{0.701}\,$\pm$\,0.012 & \textbf{0.808}\,$\pm$\,0.003 & 0.381\,$\pm$\,0.022 & \textbf{0.512}\,$\pm$\,0.007 & \textbf{16.7} & \textbf{25.1} \\
    GPT-5.6-sol & 0.679\,$\pm$\,0.035 & 0.775\,$\pm$\,0.013 & 0.362\,$\pm$\,0.007 & 0.509\,$\pm$\,0.009 & 18.3 & 32.0 \\
    Qwen3.8-27B & 0.552\,$\pm$\,0.006 & 0.694\,$\pm$\,0.010 & \textbf{0.426}\,$\pm$\,0.031 & 0.498\,$\pm$\,0.008 & 23.3 & 32.7 \\
    \bottomrule
  \end{tabular}
\end{table}

To test whether Stage~1 initialization depends on a specific model or pool size, \cref{tab:llm_init} compares $k=1$ and $k=3$ pools across three models released between July and September 2026, namely Claude Opus~5.5 and GPT-5.6-sol at high reasoning effort, and the compact open-weight Qwen3.8-27B.
We run Qwen3.8-27B with thinking disabled because thinking mode substantially increases inference latency.
At $k=3$, all three models outperform the random, expert, and hMETIS baselines in \cref{tab:main_stage_results} on both suites.
Notably, Qwen3.8-27B leads single-shot held-out proposals and closely approaches the frontier models at $k=3$, showing that effective formulation-based initialization does not require a frontier reasoning model.

Across all three models, however, individual proposals remain unreliable.
Roughly one-sixth to one-third of candidates score zero due to unbounded or intractable subproblems, and about half of Claude Opus~5.5 three-proposal pools contain at least one zero-scoring candidate.
Evaluating a small pool of three proposals rather than a single candidate absorbs these failures and substantially raises mean Dual AUC across every model and suite.

\subsection{Primal Heuristic Generation}
\label{sec:primal_recovery_results}\label{sec:solution_quality}\label{sec:lr_state_contribution}

\begin{table}
  \centering
  \caption{Primal solution quality on hard instances under a 30-second budget. Entries report mean percentage gaps over feasible calls to $U_I$ after the initial recovery stage (Init) and after improvement (Final); negative values beat $U_I$, and bold marks the lowest Final gap.}
  \label{tab:primal_recovery}
  \footnotesize
  \setlength{\tabcolsep}{4pt}
  \begin{tabular}{@{}ll rr rr rr rr@{}}
    \toprule
    & & \multicolumn{4}{c}{\textbf{LLM-generated heuristics}} & \multicolumn{4}{c}{\textbf{Generic heuristics}} \\
    \cmidrule(lr){3-6}\cmidrule(l){7-10}
    & & \multicolumn{2}{c}{\textbf{LR-aware}} & \multicolumn{2}{c}{\textbf{LR-blind}} & \multicolumn{2}{c}{\textbf{FJ+$\bar x$-RINS}} & \multicolumn{2}{c}{\textbf{$\bar x$-RINS/RENS}} \\
    \cmidrule(lr){3-4}\cmidrule(lr){5-6}\cmidrule(lr){7-8}\cmidrule(l){9-10}
    \textbf{Set} & \textbf{Bench.} & \textbf{Init} & \textbf{Final} & \textbf{Init} & \textbf{Final} & \textbf{Init} & \textbf{Final} & \textbf{Init} & \textbf{Final} \\
    \midrule
    \multirow{9}{*}{Validation}
    & P01 & 5.92 & \textbf{$-$1.90} & 1.53 & $-$1.48 & \multicolumn{2}{c}{fail} & \multicolumn{2}{c}{fail} \\
    & P02 & 38.62 & \textbf{0.85} & 23.37 & 2.06 & $>$100 & 5.25 & \multicolumn{2}{c}{fail} \\
    & P03 & 15.66 & $-$0.54 & 2.40 & \textbf{$-$0.62} & 33.47 & 30.97 & \multicolumn{2}{c}{fail} \\
    & P04 & 27.79 & 3.87 & $>$100 & 10.02$^{\dagger}$ & \multicolumn{2}{c}{fail} & 0.40 & \textbf{0.40} \\
    & P05 & 36.22 & 5.86 & 16.12 & 3.90 & $>$100 & \textbf{0.65} & \multicolumn{2}{c}{fail} \\
    & P06 & 9.46 & 1.76 & 10.98 & \textbf{1.17} & \multicolumn{2}{c}{fail} & \multicolumn{2}{c}{fail} \\
    & P07 & 8.80 & \textbf{1.09} & 21.69 & 1.23 & $>$100 & 42.76 & \multicolumn{2}{c}{fail} \\
    & P09 & 2.06 & 0.64 & 1.32 & 0.75 & 54.13 & \textbf{0.59} & 0.71 & 0.71 \\
    & P10 & 1.14 & $-$0.20 & 27.02 & 2.09 & \multicolumn{2}{c}{fail} & $-$0.52 & \textbf{$-$0.52} \\
    \cmidrule(l){2-10}
    & Mean & & \textbf{1.27} & & 2.12 & \multicolumn{2}{c}{---} & \multicolumn{2}{c}{---} \\
    \midrule
    \multirow{8}{*}{Held-out}
    & T01 & $>$100 & \textbf{27.34} & $>$100 & 41.69 & $>$100 & $>$100 & \multicolumn{2}{c}{fail} \\
    & T02 & $>$100 & 1.64 & 64.81 & \textbf{0.93} & $>$100 & $>$100 & \multicolumn{2}{c}{fail} \\
    & T04 & 10.91 & \textbf{$-$0.75} & 4.62 & \textbf{$-$0.75} & 37.90 & 37.90$^{\ddagger}$ & \multicolumn{2}{c}{fail} \\
    & T05 & 67.10 & 0.27 & 66.74 & \textbf{0.12} & 64.20 & 35.96 & 40.50 & 40.50 \\
    & T07 & 23.05 & 1.16 & 20.86 & \textbf{0.59} & \multicolumn{2}{c}{fail} & \multicolumn{2}{c}{fail} \\
    & T08 & 45.86 & \textbf{1.21} & 3.62 & 2.72 & \multicolumn{2}{c}{fail} & \multicolumn{2}{c}{fail} \\
    & T09 & $>$100 & $-$0.24 & 11.64 & \textbf{$-$1.31} & $>$100 & $>$100 & \multicolumn{2}{c}{fail} \\
    & T10 & 0.86 & 0.53 & 0.29 & \textbf{0.29} & $>$100 & $>$100 & \multicolumn{2}{c}{fail} \\
    \cmidrule(l){2-10}
    & Mean & & \textbf{3.89} & & 5.54 & \multicolumn{2}{c}{---} & \multicolumn{2}{c}{---} \\
    \bottomrule
  \end{tabular}
  \par\smallskip
  \parbox{\linewidth}{\scriptsize $^{\dagger}$8/9 feasible calls across 3 replicates $\times$ 3 instances. \quad $^{\ddagger}$2/3 feasible instances.}
\end{table}

We evaluate Phase~3 primal heuristic generation on unseen hard instances in \cref{tab:primal_recovery}, first comparing generated solver-free programs against the two solver-assisted baselines in Appendix~\ref{app:xbar_reference}, and then assessing the contribution of the LR state $s_I$ via an LR-blind ablation.
Phase~3 runs three replicates on medium instances for both the LR-aware and LR-blind agents and evaluates each selected program without modification on unseen hard instances under a single-threaded 30-second budget.
Performance is measured by percentage gap to $U_I$, excluding the preceding LR run and offline program generation, which takes about 40 minutes per replicate.
The evaluation covers the nine validation and eight held-out benchmarks whose medium-selected AutoLR relaxation achieves a positive hard-instance Dual AUC, as analyzed in \cref{sec:selected_relaxations}.

As shown in \cref{tab:primal_recovery}, the LR-aware agent recovers a feasible solution on every call across all 17 benchmarks, staying within 2\% of $U_I$ on 14 of them.
By contrast, both generic solver-assisted baselines frequently fail to restore feasibility within 30 seconds, and when they succeed, they remain competitive on a few validation benchmarks.

Comparing the Init and Final columns of \cref{tab:primal_recovery} highlights the distinct roles of initial recovery and anytime improvement.
In practice, generated $R_\phi$ functions span standard Lagrangian heuristic patterns, including iterate repair, multi-candidate selection, multiplier-guided construction, and occasionally brief local polishing. When initial gaps remain large, $Q_\psi$ drives most of the gap reduction.
For example, on the graph-drawing benchmark T02, $R_\phi$ recovers feasible layouts from LR iterates while $Q_\psi$ refines coordinates with local moves and simulated annealing, cutting the LR-aware gap from above 100\% to under 2\%.
On the GPU layer-placement benchmark T08, direct rounding of LR iterates can yield zero-flow feasible solutions, so $R_\phi$ also constructs placements guided by $s_I$ while $Q_\psi$ searches placements evaluated by maximum flow, reducing the mean gap from 45.9\% to roughly 1\%.

The LR-blind ablation shows that formulation-aware code generation drives most of the performance, whereas Lagrangian dual states help selectively on a subset of benchmarks.
Across \cref{tab:primal_recovery}, the two agents split final wins evenly with no statistically significant paired difference.
Although the LR-blind agent often reports tighter Init values because it is not constrained to depend on $s_I$ and more frequently includes brief polishing inside $R_\phi$, both agents are compared under the same 30-second total budget in Final.
Moreover, LR-blind programs can still synthesize internal Lagrangian pricing on their own.
Thus, when only a primal solution is needed, standalone code synthesis suffices without running LR, whereas in a primal--dual LR workflow $s_I$ is available as a byproduct of dual optimization and yields clear gains on benchmarks such as P04, P10, and T01 through dual prices and LR-guided priorities.

\section{Benchmark-Level Analysis}
\label{sec:instance_wise_results}
Building on the suite-level results in \cref{sec:experimental_results}, \cref{sec:selected_relaxations} examines what relaxation structures AutoLR selects relative to published expert designs, hMETIS partitions, and LP relaxations, and when medium-selected relaxations fail to scale to harder instances.
\Cref{sec:search_trajectories} then traces how Stage~2 search refines relaxations from feedback, and \cref{sec:anytime_replays} examines how primal recovery and dual bounds evolve together during an LR run.

\subsection{Selected Relaxations and Scale Transfer}
\label{sec:selected_relaxations}

\begin{table}
  \centering
  \caption{Selected relaxations on the validation benchmarks: retained structure, LR and LP bound gaps on medium instances, and Dual AUC on medium and hard instances under a 120-second LR budget.}
  \label{tab:autolr_structure_performance}
  \footnotesize
  \setlength{\tabcolsep}{4pt}
  \renewcommand{\arraystretch}{1.08}
  \begin{tabular*}{\textwidth}{@{\extracolsep{\fill}}lcl rr rr rr@{}}
    \toprule
     &  &  & \multicolumn{2}{c}{\shortstack{\textbf{Structure}\\\textbf{(medium)}}} & \multicolumn{2}{c}{\shortstack{\textbf{Bound gap (\%)}\\\textbf{(medium)}}} & \multicolumn{2}{c}{\shortstack{\textbf{Dual AUC}}} \\
    \cmidrule(lr){4-5}\cmidrule(lr){6-7}\cmidrule(lr){8-9}
    \textbf{Bench.} & \shortstack{\textbf{Total}\\\textbf{families}} & \textbf{Method} & \shortstack{\textbf{Dualized}\\\textbf{rows (\%)}} & \textbf{Blocks} & \textbf{LR} & \textbf{LP} & \textbf{Medium} & \textbf{Hard} \\
    \midrule
    \multirow{2}{*}{\textbf{P01}} & \multirow{2}{*}{5} & Expert & 0.7 & 75 & \textbf{0.2} & \multirow{2}{*}{61.3} & \textbf{0.975} & \textbf{0.951} \\
      &  & AutoLR & 0.7 & 75 & \textbf{0.2} &  & \textbf{0.975} & \textbf{0.951} \\
    \addlinespace[4pt]
    \multirow{2}{*}{\textbf{P02}} & \multirow{2}{*}{5} & Expert & 1.7 & 1 & \textbf{24.9} & \multirow{2}{*}{26.5} & \textbf{0.783} & \textbf{0.697} \\
      &  & AutoLR & 1.7 & 1 & \textbf{24.9} &  & \textbf{0.783} & \textbf{0.697} \\
    \addlinespace[4pt]
    \multirow{2}{*}{\textbf{P03}} & \multirow{2}{*}{6} & Expert & 0.8 & 1 & 22.9 & \multirow{2}{*}{11.4} & 0.733 & 0.143 \\
      &  & AutoLR & 90.2 & 1 & \textbf{11.4} &  & \textbf{0.887} & \textbf{0.889} \\
    \addlinespace[4pt]
    \multirow{2}{*}{\textbf{P04}} & \multirow{2}{*}{6} & Expert & 21.5 & 15 & 4.1 & \multirow{2}{*}{39.5} & 0.944 & \textbf{0.939} \\
      &  & AutoLR & 22.1 & 4.3 & \textbf{3.9} &  & \textbf{0.945} & 0.922 \\
    \addlinespace[4pt]
    \multirow{2}{*}{\textbf{P05}} & \multirow{2}{*}{6} & Expert & 1.1 & 1 & \textbf{5.1} & \multirow{2}{*}{12.8} & 0.820 & 0.778 \\
      &  & AutoLR & 33.9 & 1 & \textbf{5.1} &  & \textbf{0.896} & \textbf{0.871} \\
    \addlinespace[4pt]
    \multirow{2}{*}{\textbf{P06}} & \multirow{2}{*}{7} & Expert & 57.9 & 4021 & 14.8 & \multirow{2}{*}{11.4} & 0.820 & 0.522 \\
      &  & AutoLR & 38.1 & 20 & \textbf{2.6} &  & \textbf{0.942} & \textbf{0.840} \\
    \addlinespace[4pt]
    \multirow{2}{*}{\textbf{P07}} & \multirow{2}{*}{9} & Expert & 19.0 & 6 & \textbf{23.2} & \multirow{2}{*}{34.6} & 0.634 & 0.562 \\
      &  & AutoLR & 94.9 & 6 & 23.3 &  & \textbf{0.705} & \textbf{0.678} \\
    \addlinespace[4pt]
    \multirow{2}{*}{\textbf{P08}} & \multirow{2}{*}{9} & Expert & 21.7 & 5 & 70.7 & \multirow{2}{*}{81.7} & 0.245 & \textbf{0.480} \\
      &  & AutoLR & 50.9 & 1 & \textbf{57.2} &  & \textbf{0.346} & 0.000 \\
    \addlinespace[4pt]
    \multirow{2}{*}{\textbf{P09}} & \multirow{2}{*}{9} & Expert & 0.4 & 22 & \textbf{1.3} & \multirow{2}{*}{22.3} & 0.961 & 0.943 \\
      &  & AutoLR & 68.1 & 22 & \textbf{1.3} &  & \textbf{0.965} & \textbf{0.949} \\
    \addlinespace[4pt]
    \multirow{2}{*}{\textbf{P10}} & \multirow{2}{*}{10} & Expert & 0.5 & 1 & 26.0 & \multirow{2}{*}{100.0} & 0.564 & 0.000 \\
      &  & AutoLR & 5.0 & 8 & \textbf{1.9} &  & \textbf{0.961} & \textbf{0.971} \\
    \midrule
    \multirow{2}{*}{\textbf{Mean}} &  & Expert & \multicolumn{2}{c}{} &  &  & 0.748 & 0.601 \\
      &  & AutoLR & \multicolumn{2}{c}{} &  &  & \textbf{0.840} & \textbf{0.777} \\
    \bottomrule
  \end{tabular*}
  \par\smallskip
  \parbox{\linewidth}{\scriptsize Entries average three instances per level. Dualized rows is the percentage of dualized constraint rows. Bound gaps are relative to the 120\,s Gurobi incumbent. Bold marks higher Dual AUC or lower LR gap (including ties); $>$100 exceeds 100\%; --- denotes no finite bound.}
\end{table}

Selected relaxations improve Dual AUC by balancing bound strength against subproblem cost, namely solution time, rather than by maximizing block counts, and their transfer to harder instances depends on how the solution time of the retained subproblems grows with instance size.
Phase~2 of AutoLR, described in \cref{sec:autolr_selection}, consolidates the nine medium-instance LLM--LLM runs of each benchmark into one AutoLR relaxation by majority vote or highest mean Dual AUC, and reuses it unchanged on the hard instances.
\Cref{tab:autolr_structure_performance,tab:heldout_structure_performance} compare the selected AutoLR relaxations with the baselines, reporting medium-instance structure and bound gaps alongside Dual AUC at both difficulty levels.
Appendices~\ref{app:short_validation} and~\ref{app:short_heldout} summarize the dualized constraints and resulting subproblem structures for each benchmark.

\label{sec:selected_validation}
On the validation set in \cref{tab:autolr_structure_performance}, AutoLR reproduces the expert relaxation on P01 and P02 and selects different constraint families on the remaining eight benchmarks.
Seven of the eight differing selections still dualize at least one family from the expert relaxation.
Across the ten benchmarks, AutoLR raises medium Dual AUC on seven, matches or practically ties the expert relaxation on the other three, and tightens or preserves the medium LR bound gap on nine.
These gains over the expert relaxations, however, reflect finite-budget adaptation to a general-purpose LR engine rather than superiority over the source papers' specialized solvers, as expert relaxations are often tailored to dedicated subproblem algorithms.
Independently of this comparison, AutoLR improves upon the LP bound on nine of the ten benchmarks.

The largest Dual AUC gains occur where the expert relaxation hits a computational bottleneck within the 120-second budget, either in costly subproblems or in bundle overhead.
On P03 and P10, the expert relaxation leaves a single coupled subproblem block, typically intended for specialized combinatorial or dynamic-programming algorithms, that stalls in the generic MILP solver.
On P03, that block slows and times out after a few iterations, whereas AutoLR dualizes most rows to lower its cost and halves the medium bound gap.
On P10, the expert block is costly from the first iteration, whereas AutoLR dualizes machine capacity to expose the per-job subproblems in \cref{fig:selected_structures_p10}, cutting a 100\% LP gap to under 2\%.
On P06, the opposite bottleneck occurs because fine splitting pushes most of the budget into bundle overhead and stalls the expert above the LP gap.
AutoLR selects the coarser decomposition in \cref{fig:selected_structures_p06}, shifting runtime back to subproblem evaluations and tightening the gap well below the LP gap.

AutoLR also raises Dual AUC by lowering subproblem cost when dualizing additional families preserves both the block count and bound strength.
On P05, P07, and P09, AutoLR dualizes a superset of the expert families, retains the same number of blocks, and reaches the same medium bound gap within 0.1 percentage points.
By \cref{obs:monotone}, its optimal dual value cannot exceed the expert bound, so its higher Dual AUC comes from faster subproblem evaluations that advance the dual bound earlier within the budget.
Because Dual AUC integrates the best-so-far bound over time without directly penalizing dualized rows, extra multipliers matter only through bundle overhead.
Among candidates reaching the same bound, Dual AUC therefore favors the one with the lowest subproblem cost, even when it dualizes far more rows than the expert relaxation.

Where the expert relaxation already solves cleanly in multiple blocks, AutoLR can tighten medium bounds by retaining larger blocks, as on P04 and P08.
On P04, AutoLR dualizes a different set of families that leaves fewer and larger blocks, marginally tightening the medium bound and practically tying medium Dual AUC.
On P08, it similarly leaves a single coupled block instead of five, substantially tightening the medium bound gap and raising medium Dual AUC.

\label{sec:selected_heldout}
\begin{table}
  \centering
  \caption{Selected relaxations on the held-out benchmarks: retained structure, LR and LP bound gaps on medium instances, and Dual AUC on medium and hard instances under a 120-second LR budget.}
  \label{tab:heldout_structure_performance}
  \footnotesize
  \setlength{\tabcolsep}{4pt}
  \renewcommand{\arraystretch}{1.08}
  \begin{tabular*}{\textwidth}{@{\extracolsep{\fill}}lcl rr rr rr@{}}
    \toprule
     &  &  & \multicolumn{2}{c}{\shortstack{\textbf{Structure}\\\textbf{(medium)}}} & \multicolumn{2}{c}{\shortstack{\textbf{Bound gap (\%)}\\\textbf{(medium)}}} & \multicolumn{2}{c}{\shortstack{\textbf{Dual AUC}}} \\
    \cmidrule(lr){4-5}\cmidrule(lr){6-7}\cmidrule(lr){8-9}
    \textbf{Bench.} & \shortstack{\textbf{Total}\\\textbf{families}} & \textbf{Method} & \shortstack{\textbf{Dualized}\\\textbf{rows (\%)}} & \textbf{Blocks} & \textbf{LR} & \textbf{LP} & \textbf{Medium} & \textbf{Hard} \\
    \midrule
    \multirow{2}{*}{\textbf{T01}} & \multirow{2}{*}{7} & hMETIS & 0.3 & 2.3 & \textbf{21.1} & \multirow{2}{*}{$>$100} & \textbf{0.802} & \textbf{0.841} \\
      &  & AutoLR & 68.2 & 1 & 43.1 &  & 0.685 & 0.733 \\
    \addlinespace[4pt]
    \multirow{2}{*}{\textbf{T02}} & \multirow{2}{*}{15} & hMETIS & 7.0 & 2 & 94.1 & \multirow{2}{*}{15.8} & 0.059 & 0.056 \\
      &  & AutoLR & 38.5 & 1 & \textbf{14.0} &  & \textbf{0.852} & \textbf{0.827} \\
    \addlinespace[4pt]
    \multirow{2}{*}{\textbf{T03}} & \multirow{2}{*}{7} & hMETIS & 32.0 & 13 & $>$100 & \multirow{2}{*}{$>$100} & 0.000 & 0.000 \\
      &  & AutoLR & --- & --- & --- &  & --- & --- \\
    \addlinespace[4pt]
    \multirow{2}{*}{\textbf{T04}} & \multirow{2}{*}{5} & hMETIS & 88.2 & 11 & 43.3 & \multirow{2}{*}{43.3} & 0.542 & 0.594 \\
      &  & AutoLR & 29.4 & 11 & \textbf{32.1} &  & \textbf{0.671} & \textbf{0.735} \\
    \addlinespace[4pt]
    \multirow{2}{*}{\textbf{T05}} & \multirow{2}{*}{7} & hMETIS & 0.3 & 49 & 10.4 & \multirow{2}{*}{10.4} & 0.879 & \textbf{0.908} \\
      &  & AutoLR & 16.3 & 1 & \textbf{5.3} &  & \textbf{0.928} & 0.740 \\
    \addlinespace[4pt]
    \multirow{2}{*}{\textbf{T06}} & \multirow{2}{*}{19} & hMETIS & 9.5 & 2 & 18.8 & \multirow{2}{*}{100.0} & 0.799 & 0.000 \\
      &  & AutoLR & 38.6 & 1 & \textbf{1.6} &  & \textbf{0.959} & 0.000 \\
    \addlinespace[4pt]
    \multirow{2}{*}{\textbf{T07}} & \multirow{2}{*}{19} & hMETIS & 3.7 & 74 & \textbf{1.3} & \multirow{2}{*}{42.7} & \textbf{0.832} & 0.000 \\
      &  & AutoLR & 30.7 & 114.3 & 36.4 &  & 0.612 & \textbf{0.463} \\
    \addlinespace[4pt]
    \multirow{2}{*}{\textbf{T08}} & \multirow{2}{*}{10} & hMETIS & 24.6 & 2 & $>$100 & \multirow{2}{*}{$>$100} & 0.022 & 0.031 \\
      &  & AutoLR & 37.5 & 1 & $>$100 &  & \textbf{0.046} & \textbf{0.065} \\
    \addlinespace[4pt]
    \multirow{2}{*}{\textbf{T09}} & \multirow{2}{*}{22} & hMETIS & 11.2 & 2 & --- & \multirow{2}{*}{93.5} & 0.000 & 0.000 \\
      &  & AutoLR & 50.4 & 1 & \textbf{93.5} &  & \textbf{0.064} & \textbf{0.157} \\
    \addlinespace[4pt]
    \multirow{2}{*}{\textbf{T10}} & \multirow{2}{*}{10} & hMETIS & 40.1 & 3 & 57.0 & \multirow{2}{*}{23.7} & 0.423 & 0.782 \\
      &  & AutoLR & 55.1 & 4 & \textbf{10.3} &  & \textbf{0.885} & \textbf{0.905} \\
    \midrule
    \multirow{2}{*}{\textbf{Mean}} &  & hMETIS & \multicolumn{2}{c}{} &  &  & 0.436 & 0.321 \\
      &  & AutoLR & \multicolumn{2}{c}{} &  &  & \textbf{0.570} & \textbf{0.463} \\
    \bottomrule
  \end{tabular*}
  \par\smallskip
  \parbox{\linewidth}{\scriptsize Definitions and formatting follow \cref{tab:autolr_structure_performance}. hMETIS may dualize individual rows within a family, with structure reported regardless of bound finiteness. On T03, AutoLR selects no relaxation (---), counting as zero in Dual AUC means.}
\end{table}

On the held-out set in \cref{tab:heldout_structure_performance}, AutoLR outperforms hMETIS on most benchmarks because a purely structural partition ignores the algebraic role of the relaxed constraints.
AutoLR achieves higher medium Dual AUC than hMETIS on seven of the ten benchmarks and trails on T01 and T07.
On T03, AutoLR selects no relaxation because every medium candidate scores zero, as does hMETIS.
Moreover, AutoLR improves upon the LP bound on seven of the ten benchmarks without relying on published LR designs.
Unlike hMETIS, which splits every benchmark into multiple blocks while minimizing dualized rows, AutoLR dualizes a larger share of rows on eight of its nine selected relaxations and retains a single block on six.

Where hMETIS dualizes fewer rows and forms more blocks, AutoLR still achieves higher Dual AUC and tighter bounds with a single retained block, as on T02, T05, and T06.
On T02, hMETIS relaxes a bounding constraint that soon leaves the subproblem unbounded, whereas the single-block AutoLR relaxation in \cref{fig:selected_structures_t02} stays bounded and improves upon the LP bound.
On T05, hMETIS relaxes the linking families to form many small blocks but collapses to the LP bound, whereas AutoLR keeps one coupled block that nearly halves the LP gap.
On T06, splitting coupled decisions leaves hMETIS with a wide bound gap, whereas AutoLR retains a single block close to the original model that solves in seconds and nearly closes the LP gap.

Even when AutoLR and hMETIS form a similar number of blocks, AutoLR obtains stronger bounds and lower subproblem cost by dualizing more suitable constraints.
On T04, both methods retain 11 blocks, but hMETIS dualizes most rows across both linking families and falls back to the LP bound, whereas AutoLR relaxes only one linking family, preserving pairwise coupling within each block to improve upon the LP bound.
On T10, AutoLR relaxes feasibility-only constraints so the subproblem retaining the objective terms solves within seconds and improves upon the LP bound.
By contrast, hMETIS also relaxes rows coupled to the objective and leaves costly subproblems that hit the time limit on some instances, keeping its bound gap wider than the LP gap.

Row-level partitioning lets hMETIS outperform AutoLR on medium instances of T01 and T07, outside the family-level design space.
On T01, AutoLR identifies the strongest family-level relaxation by dualizing a single large constraint family into one retained block, whereas hMETIS dualizes less than 1\% of the rows across subsets of three families, forms multiple blocks, and nearly halves the medium bound gap.
On T07, hMETIS similarly leads by relaxing a few rows that cut the schedule into multi-day spans.

Some formulations resist Lagrangian relaxation regardless of the selection method, leaving both AutoLR and hMETIS below 0.1 medium Dual AUC.
On T03, big-M rows supply the only upper bounds on some variables, so relaxing them yields unbounded subproblems whereas keeping them leaves negative bounds that score zero in Dual AUC.
On T08, strong bounds require an intact chain of constraint families, so any relaxation with low subproblem cost breaks that chain and caps the bound far below what Gurobi certifies within seconds.
On T09, shared variables couple the blocks under every family-level mask, so making the benchmark amenable to LR requires a reformulation that splits those variables by block.

\label{sec:selected_transfer}
Applied unchanged to the hard instances, the medium-selected relaxations keep most of their advantage, while a few retained structures fail to scale.
On the validation set, AutoLR matches or exceeds the expert relaxation on hard instances of eight of the ten benchmarks and raises mean hard Dual AUC by 0.176, although the ten-benchmark paired difference does not reach the 5\% level.
The gain concentrates on P03 and P10, where the single expert block stalls on the hard instances and keeps little or none of its medium score.
On the held-out set, AutoLR raises mean hard Dual AUC over hMETIS by 0.142, winning on six benchmarks, trailing on T01 and T05, and tying at zero on T03 and T06.
Most of these differences stem from whether a relaxation obtains a finite bound and how tightly it closes the gap.

Transfer failures on hard instances trace to how the solution time of the retained subproblems grows with instance size.
On P08 and T06, AutoLR's single-block selections exceed 120 seconds on the first hard subproblem and drop to zero.
On T05, AutoLR's larger blocks still yield finite hard bounds, but slower iterations and delayed convergence pull AutoLR below the multi-block baseline.
Even when recomputed on each hard instance, hMETIS suffers the same failure on T07, where its blocks exceed 120 seconds while AutoLR's finer decomposition remains tractable.

Selection on medium instances tends to favor larger blocks whenever they solve quickly at that scale, yet the search itself still surfaces scalable designs.
Checking the top candidates at the target scale before fixing a single relaxation can therefore avoid some of these failures.
Selecting relaxations directly from hard-instance runs restores strong hard performance on P08, T05, and T06 and finds positive-scoring masks on hard T03.
Although those masks are selected on the same hard instances that score them, the scalable alternatives for P08 and T05 already appear during medium-instance search, trailing the medium winners only slightly.

\subsection{Feedback-Driven Search Trajectories}
\label{sec:search_trajectories}

\begin{figure}
  \centering
  \captionsetup{font=small}
  \captionsetup[subfigure]{skip=1pt}
  \makebox[0.195\textwidth][c]{\small\textbf{Expert}}\hfill
  \makebox[0.195\textwidth][c]{\small\textbf{AutoLR}}\hfill
  \makebox[0.195\textwidth][c]{\small\textbf{hMETIS}}\hfill
  \makebox[0.195\textwidth][c]{\small\textbf{Search}}\par
  \begin{subfigure}[t]{\textwidth}
    \centering
    \begin{minipage}[t]{0.195\textwidth}\centering
      \includegraphics[width=\linewidth]{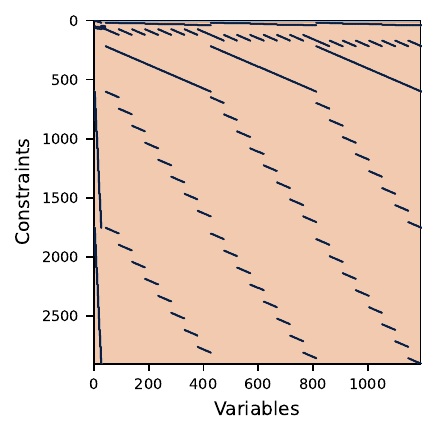}\par
      {\scriptsize 16 rows, 1 block\par}
    \end{minipage}\hfill
    \begin{minipage}[t]{0.195\textwidth}\centering
      \includegraphics[width=\linewidth]{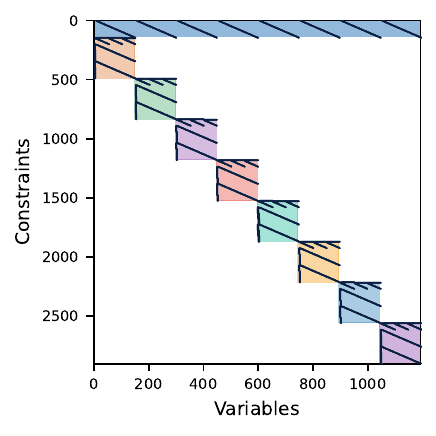}\par
      {\scriptsize 144 rows, 8 blocks\par}
    \end{minipage}\hfill
    \begin{minipage}[t]{0.195\textwidth}\centering
      \includegraphics[width=\linewidth]{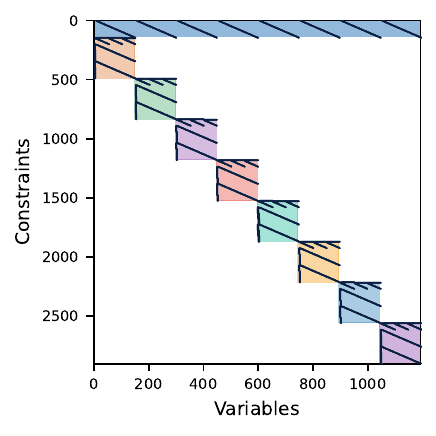}\par
      {\scriptsize Same as AutoLR\par}
    \end{minipage}\hfill
    \begin{minipage}[t]{0.195\textwidth}\centering
      \includegraphics[width=\linewidth]{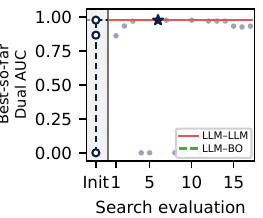}\par
      {\scriptsize Pool / final: 0.976 / 0.976\par}
    \end{minipage}
    \caption{P10: the selected relaxation is already in the initial pool.}
    \label{fig:selected_structures_p10}
  \end{subfigure}\\[0pt]
  \begin{subfigure}[t]{\textwidth}
    \centering
    \begin{minipage}[t]{0.195\textwidth}\centering
      \includegraphics[width=\linewidth]{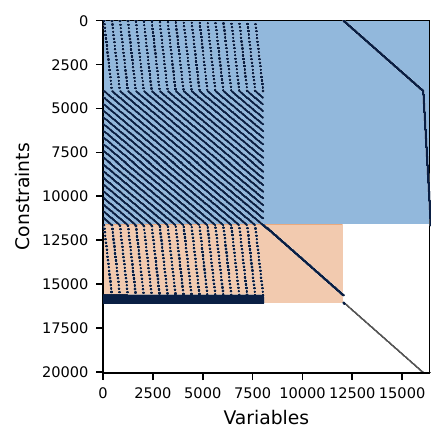}\par
      {\scriptsize 11,600 rows, 4,021 blocks\par}
    \end{minipage}\hfill
    \begin{minipage}[t]{0.195\textwidth}\centering
      \includegraphics[width=\linewidth]{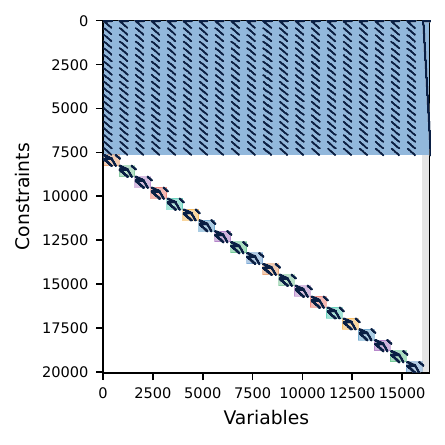}\par
      {\scriptsize 7,640 rows, 20 blocks\par}
    \end{minipage}\hfill
    \begin{minipage}[t]{0.195\textwidth}\centering
      \includegraphics[width=\linewidth]{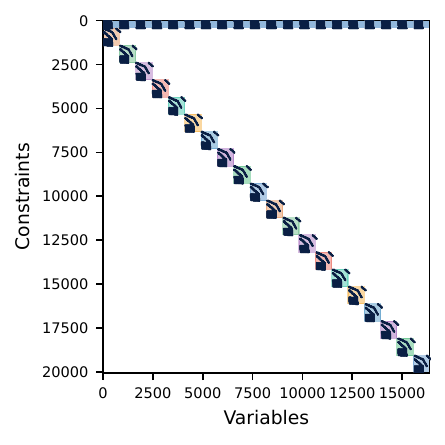}\par
      {\scriptsize 440 rows, 20 blocks\par}
    \end{minipage}\hfill
    \begin{minipage}[t]{0.195\textwidth}\centering
      \includegraphics[width=\linewidth]{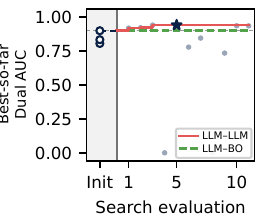}\par
      {\scriptsize Pool / final: 0.898 / 0.940\par}
    \end{minipage}
    \caption{P06: search finds the selected relaxation at evaluation 5.}
    \label{fig:selected_structures_p06}
  \end{subfigure}\\[0pt]
  \begin{subfigure}[t]{\textwidth}
    \centering
    \begin{minipage}[t]{0.195\textwidth}\centering
      \parbox[b][0.21\textwidth][c]{\linewidth}{\centering\scriptsize\itshape No published relaxation}\par
    \end{minipage}\hfill
    \begin{minipage}[t]{0.195\textwidth}\centering
      \includegraphics[width=\linewidth]{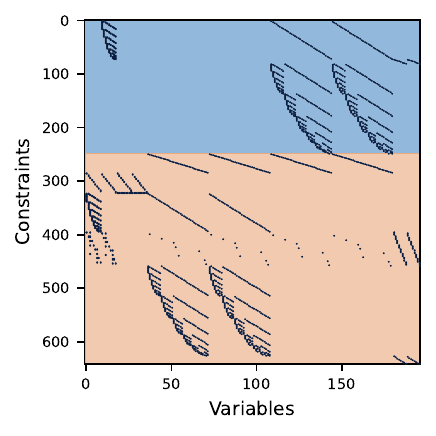}\par
      {\scriptsize 248 rows, 1 block\par}
    \end{minipage}\hfill
    \begin{minipage}[t]{0.195\textwidth}\centering
      \includegraphics[width=\linewidth]{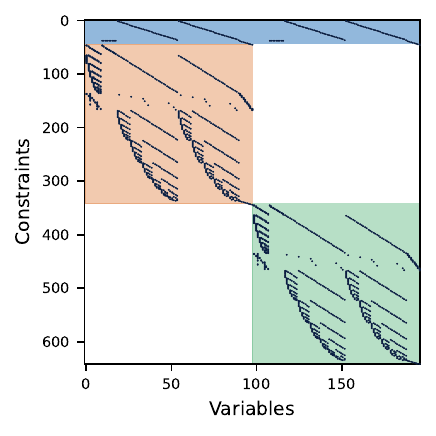}\par
      {\scriptsize 45 rows, 2 blocks\par}
    \end{minipage}\hfill
    \begin{minipage}[t]{0.195\textwidth}\centering
      \includegraphics[width=\linewidth]{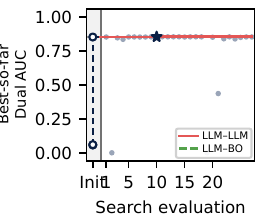}\par
      {\scriptsize Pool / final: 0.852 / 0.853\par}
    \end{minipage}
    \caption{T02: search provides a small improvement over the initial pool.}
    \label{fig:selected_structures_t02}
  \end{subfigure}\\[0pt]
  \begin{subfigure}[t]{\textwidth}
    \centering
    \begin{minipage}[t]{0.195\textwidth}\centering
      \parbox[b][0.21\textwidth][c]{\linewidth}{\centering\scriptsize\itshape No published relaxation}\par
    \end{minipage}\hfill
    \begin{minipage}[t]{0.195\textwidth}\centering
      \includegraphics[width=\linewidth]{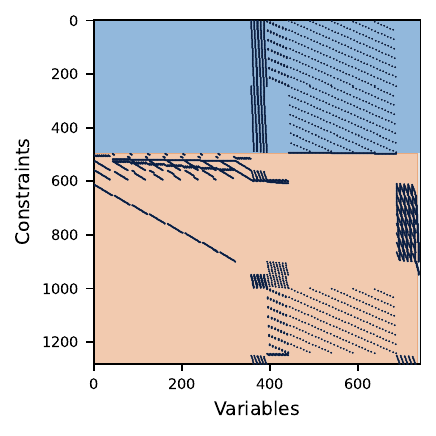}\par
      {\scriptsize 495 rows, 1 block\par}
    \end{minipage}\hfill
    \begin{minipage}[t]{0.195\textwidth}\centering
      \includegraphics[width=\linewidth]{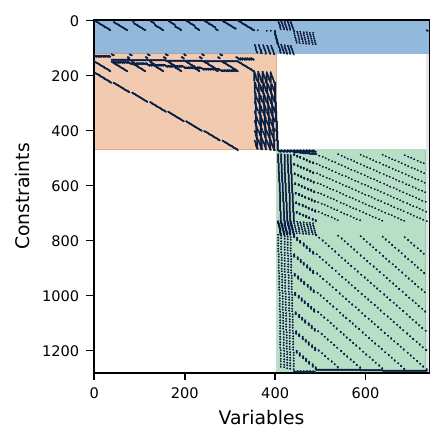}\par
      {\scriptsize 122 rows, 2 blocks\par}
    \end{minipage}\hfill
    \begin{minipage}[t]{0.195\textwidth}\centering
      \includegraphics[width=\linewidth]{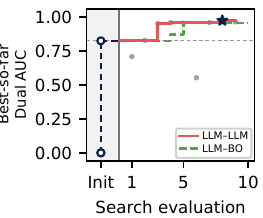}\par
      {\scriptsize Pool / final: 0.824 / 0.973\par}
    \end{minipage}
    \caption{T06: search finds the selected relaxation at evaluation 8.}
    \label{fig:selected_structures_t06}
  \end{subfigure}
  \caption{Decomposition structures and search trajectories on matched medium instances, showing the first medium run for P10 and T02 and the largest search-gain run for P06 and T06. Blue rows mark dualized constraints and colored blocks mark retained subproblems. Circles and dots mark initial and Stage~2 evaluations, solid and dashed curves show best-found Dual AUC for LLM--LLM and LLM--BO, and stars identify the displayed AutoLR relaxations.}
  \label{fig:structures_with_search}
\end{figure}

\begin{figure}
  \centering
  \includegraphics[width=\textwidth]{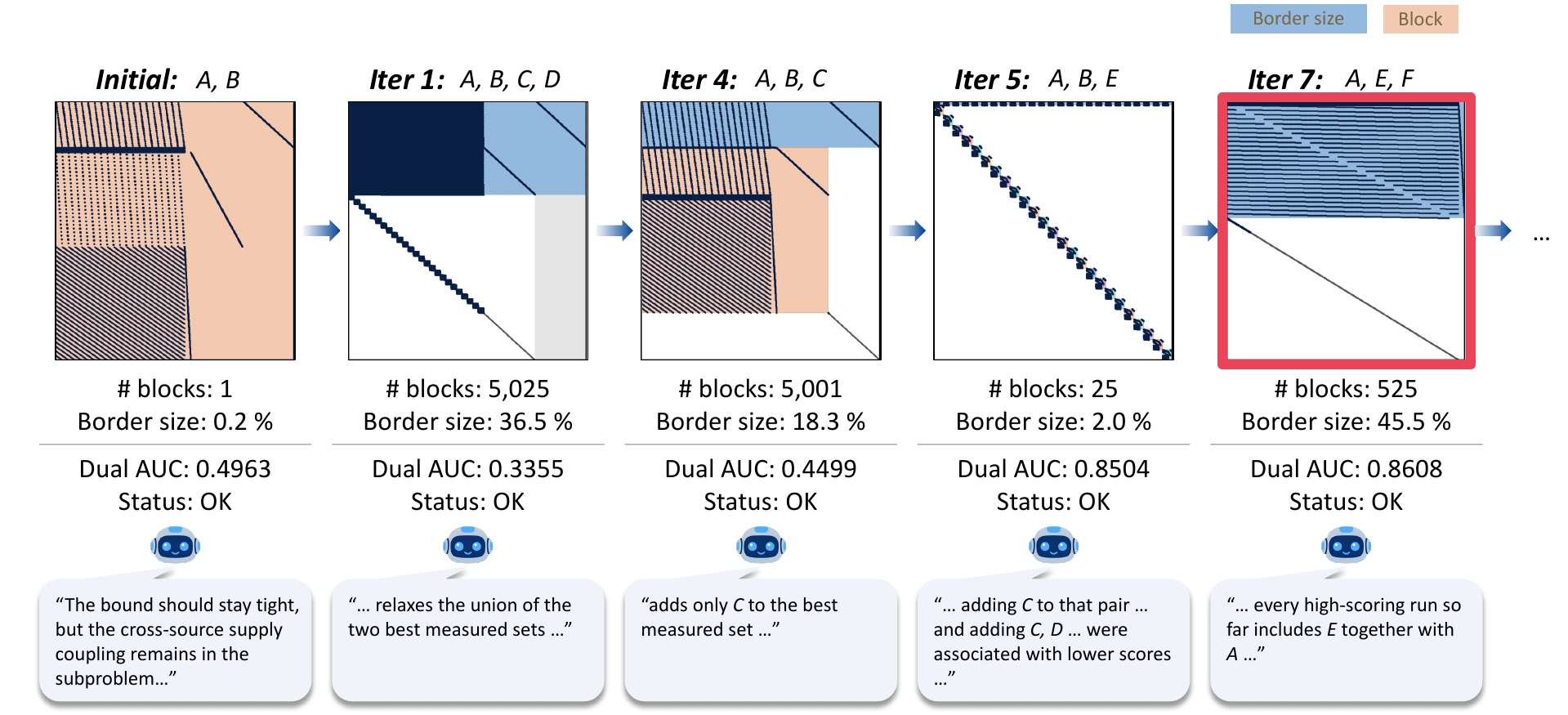}
  \caption{Feedback-driven search on a hard instance of P06 across the best initial relaxation and Stage~2 iterations 1, 4, 5, and 7, with LLM rationale excerpts. Letters A--F denote distinct relaxed constraint families, and the red frame highlights the best relaxation.}
  \label{fig:feedback_timeline}
\end{figure}

This subsection examines how Stage~2 search uses computational feedback to refine the initial pool.
\Cref{fig:structures_with_search} compares decomposition structures and search curves across four medium instances, while \cref{fig:feedback_timeline} traces one hard-instance trajectory step by step with LLM rationales.
Across the matched medium instances in \cref{fig:structures_with_search}, search confirms an already strong initial pool on P10 and T02 while achieving substantial search gains on P06 and T06.
On hard P06 in \cref{fig:feedback_timeline}, the first four iterations fail to improve upon the single-block initial relaxation.
At iteration~5, adding the single-destination family~E splits the subproblem into 25 per-store blocks and raises Dual AUC from below 0.50 to 0.85, which iteration~7 further refines to 0.86.

\subsection{Joint Primal--Dual Trajectories}
\label{sec:anytime_replays}

\begin{figure}
  \centering
  \captionsetup[subfigure]{font=small}
  \begin{subfigure}[t]{0.49\textwidth}
    \centering
    \includegraphics[width=\linewidth]{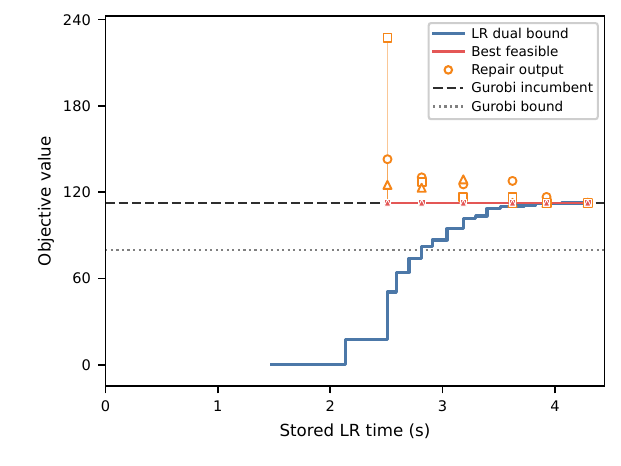}
    \caption{P10 (minimization)}
    \label{fig:anytime_p10}
  \end{subfigure}\hfill
  \begin{subfigure}[t]{0.49\textwidth}
    \centering
    \includegraphics[width=\linewidth]{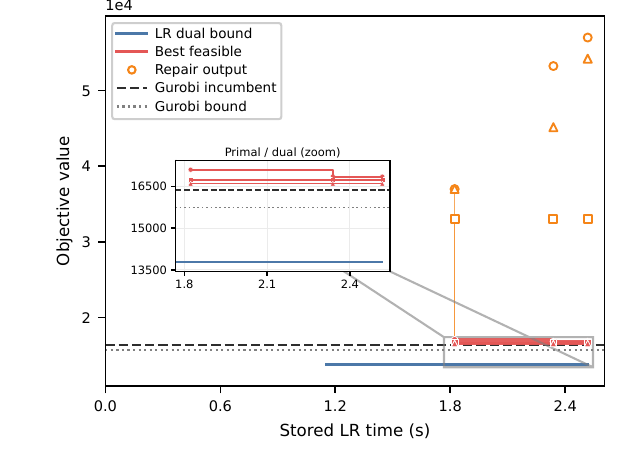}
    \caption{T02 (minimization)}
    \label{fig:anytime_t02}
  \end{subfigure}
  \caption{Joint primal--dual trajectories on hard instances of P10 and T02 across three Phase~3 replicates evaluated at intermediate LR checkpoints, excluding heuristic runtime on the horizontal axis.}
  \label{fig:anytime}
\end{figure}

To examine how primal recovery and dual bounds evolve together during an LR run rather than only at its terminal state, \cref{fig:anytime} evaluates the generated heuristics at intermediate checkpoints on unseen hard instances of P10 and T02.
On P10, as the LR run progresses over its first few seconds, $R_\phi$ outputs tighten directly onto the Gurobi incumbent while the rising LR dual bound surpasses the 120-second Gurobi bound and nearly closes the primal--dual gap.
On T02, local improvement $Q_\psi$ overcomes large $R_\phi$ gaps to approach the Gurobi incumbent from the earliest checkpoint, yet the LR dual bound remains well below those feasible values.
Together, these examples illustrate that the generated heuristics can be effective well before the full 120-second LR budget, either closing the primal--dual gap rapidly alongside a tight relaxation on P10 or delivering near-incumbent primal solutions early on T02 despite a loose dual bound.

\section{Discussion and Conclusion}
\label{sec:discussion}\label{sec:conclusion}
This paper introduced AutoLR, an agentic workflow that automates constraint-family Lagrangian relaxation design and solver-free primal heuristic generation from algebraic MILP formulations.
Across the evaluated benchmarks, LLMs add the greatest value where formulation understanding is essential. Paired with solver evaluation and feasibility checks, reading the algebraic model and synthesizing problem-specific code enable LLMs to propose strong initial relaxations and effective primal heuristics.
Once a problem-specific design is initialized, classical search methods such as Bayesian optimization remain competitive with LLM-guided refinement, and Lagrangian dual information provides selective, problem-dependent gains in primal recovery.
Broadly, these results suggest that LLMs can help reduce the manual engineering effort required to apply decomposition methods to new formulations.

Evaluating agentic decomposition designs under finite computational budgets also requires metrics aligned with anytime solver progress.
By integrating normalized dual-bound progress over time rather than scoring static structure or asymptotic bounds, Dual AUC naturally rewards relaxations that tighten bounds quickly within the budget and produces rankings that are less sensitive to the chosen time cutoff than endpoint bounds.
However, because Dual AUC measures bound quality relative to a reference incumbent rather than improvement over the LP relaxation, a high score does not by itself establish that LR substantially strengthens the LP bound.

From a deployment perspective, designing both the relaxation mask and the primal heuristic at the formulation level allows offline LLM costs to be amortized across unseen instances with zero online LLM calls.
Most medium-selected relaxations and generated heuristics transfer effectively to harder instances, although subproblem solution times grow sharply with block size, so a few coupled blocks stall at larger scales.
Even in those cases, relaxations that remain tractable on larger instances often appear among the top candidates found on medium instances. In practice, one can therefore run search on smaller instances at low cost and briefly test only the top few candidates on larger instances before fixing the final design.

\label{sec:limitations}
This study has two main limitations.
First, AutoLR selects which constraint families to dualize within a fixed formulation rather than reformulating the model itself, such as splitting shared variables across blocks or strengthening variable bounds.
Second, the workflow holds the generic subproblem and bundle solvers fixed and runs relaxation selection before primal recovery, without co-designing specialized subproblem algorithms or feeding primal and scaling outcomes back into relaxation selection.

Looking ahead, even as LLMs become faster and less expensive, testing candidate relaxations and scaling subproblems to large instances will still be limited by solver runtime rather than LLM cost.
Future work can extend AutoLR to automatic model reformulation, specialized subproblem algorithms, and feedback between primal recovery and relaxation selection, as well as related methods such as column generation and Benders decomposition.

\section*{Acknowledgments}
This work was supported by the National Research Foundation of Korea (NRF) grant funded by the Korean government (MSIT) (RS-2023-00259550).

\bibliographystyle{apalike}
\bibliography{reference}

\clearpage
\appendix
\counterwithin{table}{section}
\counterwithin{figure}{section}
\counterwithin{equation}{section}

\section{AutoLR Procedure}
\label{app:autolr_alg}

\begin{algorithm}[H]
\caption{AutoLR procedure for relaxation selection and primal heuristic generation. Dual AUC $S$, structural preview $\pi$, and mean primal--dual gap $\bar\gamma$ are defined in \cref{sec:autolr_evaluator,sec:design_and_search,sec:primal_recovery_agent}.}
\label{alg:autolr}
\small
\begin{algorithmic}[1]
\REQUIRE formulation $F$, instances $\mathcal I$, LR time budget $T$, evaluation cap $N$, replicates $R$
\ENSURE selected relaxation $z^\ast$ and primal heuristic $P^\ast$
\STATE \uline{\textbf{Phase 1: relaxation search on each instance in $\mathcal I$}}
\FORALL{$I\in\mathcal I$ and $r=1,\ldots,R$}
  \STATE $\mathcal C\leftarrow\textsc{ProposeInitial}(F)$ \COMMENT{three distinct masks fixed before evaluation}
  \STATE $\mathcal H\leftarrow\{(z,S(z;I,T),\pi(z;I)):z\in\mathcal C\}$ \COMMENT{evaluated masks with scores and previews}
  \STATE $M\leftarrow\textsc{FormulationMap}(F)$ \COMMENT{semantic map of objective and constraint families}
  \WHILE{$|\mathcal H|<N$ and the stopping rule is not met}
    \STATE $\mathcal C\leftarrow\textsc{Propose}(M,\mathcal H)$ \COMMENT{shortlist of unevaluated masks}
    \STATE $z\leftarrow\textsc{Select}(\{(c,\pi(c;I)):c\in\mathcal C\},\mathcal H)$ \COMMENT{proposal rationales withheld}
    \STATE $\mathcal H\leftarrow\mathcal H\cup\{(z,S(z;I,T),\pi(z;I))\}$
  \ENDWHILE
  \STATE $z_{I,r}\leftarrow\arg\max_{(z,s,\pi)\in\mathcal H} s$
\ENDFOR
\vspace{0.3em}
\STATE \uline{\textbf{Phase 2: relaxation selection for formulation $F$}}
\IF{$\max_{I\in\mathcal I,\,1\le r\le R} S(z_{I,r};I,T)=0$}
  \RETURN $(\varnothing,\varnothing)$ \COMMENT{no positive Dual AUC on $\mathcal I$}
\ELSIF{some mask $z$ satisfies $|\{(I,r):z_{I,r}=z\}|>|\mathcal I|R/2$}
  \STATE $z^\ast\leftarrow z$
\ELSE
  \STATE $\mathcal Z_{\mathrm{ret}}\leftarrow\{z_{I,r}:I\in\mathcal I,\ r=1,\ldots,R\}$
  \STATE $z^\ast\leftarrow\arg\max_{z\in\mathcal Z_{\mathrm{ret}}}\frac{1}{|\mathcal I|}\sum_{I\in\mathcal I} S(z;I,T)$ \COMMENT{evaluate missing $(z,I)$ pairs and reuse cached scores}
\ENDIF
\vspace{0.3em}
\STATE \uline{\textbf{Phase 3: primal heuristic generation for $z^\ast$}}
\STATE run LR with $z^\ast$ under budget $T$ on each $I\in\mathcal I$ to record $(s_I,d_I)$
\STATE $(\sigma_1,\sigma_2,\sigma_3)\leftarrow\textsc{Plan}(F,z^\ast,\{s_I\}_{I\in\mathcal I})$ \COMMENT{three distinct heuristic strategies}
\STATE $P_k\leftarrow\textsc{Implement}(F,z^\ast,\{s_I\}_{I\in\mathcal I},\sigma_k)$ for $k=1,2,3$
\STATE $\mathcal E\leftarrow\varnothing$, $P^\ast\leftarrow\varnothing$ \COMMENT{checked programs and evaluation records}
\FOR{$k=1,2,3$}
  \STATE $(P,e)\leftarrow\textsc{Check}(P_k,\{s_I,d_I\}_{I\in\mathcal I})$, $\mathcal E\leftarrow\mathcal E\cup\{(P,e)\}$
  \IF{$e$ passes and ($P^\ast=\varnothing$ or $\bar\gamma(P)<\bar\gamma(P^\ast)$)}
    \STATE $P^\ast\leftarrow P$
  \ENDIF
\ENDFOR
\IF{$P^\ast=\varnothing$}
  \RETURN $(z^\ast,\varnothing)$ \COMMENT{no initial program passes}
\ENDIF
\FOR{$k=1,2,3$}
  \STATE $P'\leftarrow\textsc{Revise}(F,z^\ast,\{s_I\}_{I\in\mathcal I},P^\ast,\mathcal E)$
  \STATE $(P,e)\leftarrow\textsc{Check}(P',\{s_I,d_I\}_{I\in\mathcal I})$, $\mathcal E\leftarrow\mathcal E\cup\{(P,e)\}$
  \IF{$e$ passes and $\bar\gamma(P)<\bar\gamma(P^\ast)$}
    \STATE $P^\ast\leftarrow P$
  \ENDIF
\ENDFOR
\RETURN $(z^\ast,P^\ast)$
\end{algorithmic}
\end{algorithm}

\clearpage
\section{Baseline Implementations}
\label{app:baselines}

\subsection{Bayesian Optimization for Relaxation Search}
\label{app:short_bo}
The BO baseline searches the set $\mathcal Z$ of nonempty masks with at most $k_{\max}$ dualized families using a zero-mean Gaussian-process (GP) surrogate of the Dual AUC $S(z;I,T)$ with the Tanimoto kernel $k_{\mathrm{T}}(z,z')=z^\top z'/(z^\top z+z'^\top z'-z^\top z')$, as in GP optimization over molecular fingerprints \citep{griffiths2023gauche}.
This kernel transfers information between overlapping masks without using family sizes or incidence structure.
At round $t$, given $n_t$ evaluated pairs $(z^{(i)},S(z^{(i)};I,T))$, $i=1,\dots,n_t$, the GP uses covariance $K=\sigma_f^2K_{\mathrm{T}}+\sigma_n^2\mathbf I_{n_t}$ with $[K_{\mathrm{T}}]_{ij}=k_{\mathrm{T}}(z^{(i)},z^{(j)})$, where $\sigma_f^2$ and $\sigma_n^2$ maximize the log marginal likelihood.
The next mask $z^{(n_t+1)}$ is chosen by enumerating $\mathcal Z\setminus\{z^{(i)}\}_{i=1}^{n_t}$ to maximize the upper confidence bound $\mu_t(z)+\sqrt{\beta}\,\sigma_t(z)$ with $\beta=2$, where $\mu_t$ and $\sigma_t$ are the GP posterior mean and standard deviation.

\subsection{Solver-Assisted Primal Baselines}
\label{app:xbar_reference}
Both solver-assisted baselines adapt RENS \citep{Berthold2014RENS} and RINS \citep{Danna2005RINS} using the bundle aggregate primal point $\bar x_\tau=\sum_{j\in\mathcal B_\tau}\nu_j x_j^{\mathrm{LR}}$ in place of an LP solution, under the same LR state and 30-second budget.
To prevent falling back to solving the unrestricted MILP, both enforce SCIP's default fixing-rate thresholds under a $10^{-5}$ integrality tolerance.
In initial recovery, \emph{$\bar x$-RINS/RENS} fixes integral discrete variables in $\bar x_\tau$, bounds the rest to $[\lfloor\bar x_{\tau,i}\rfloor,\lceil\bar x_{\tau,i}\rceil]$, and runs RENS only if at least 50\% are fixed, with one fallback attempt releasing coordinates where $\bar x_\tau$ and $x_\tau^{\mathrm{LR}}$ disagree if infeasible, whereas \emph{FJ+$\bar x$-RINS} runs Feasibility Jump \citep{Luteberget2023FeasibilityJump} from $x_\tau^{\mathrm{LR}}$ until its first feasible point.
Both then run RINS improvement from a feasible $\widehat x$ only if $\bar x_\tau$ and $\widehat x$ agree on at least 30\% of discrete variables, retaining $\widehat x$ otherwise.
Each Gurobi sub-MIP enforces a 1\% improvement cutoff on $\min\{\mathrm{UB}-\mathrm{LB},|\mathrm{UB}|\}$ within the remaining budget, and all outputs pass the verifier in \cref{sec:primal_recovery_agent}.

\section{LLM Call Structure, Token Usage, and Cost}
\label{app:short_llm_cost}

\begingroup
\setlength{\intextsep}{6pt}
\begin{table}[H]
  \centering
  \caption{Mean LLM usage per run of both agents.}
  \label{tab:short_llm_usage}
  \footnotesize
  \setlength{\tabcolsep}{2pt}
  \begin{tabular*}{\textwidth}{@{\extracolsep{\fill}}lrrrrrrrr@{}}
    \toprule
    & \multicolumn{4}{c}{\textbf{Relaxation Design Agent}} & \multicolumn{4}{c}{\textbf{Primal Heuristic Agent}} \\
    \cmidrule(lr){2-5}\cmidrule(l){6-9}
    & \multicolumn{2}{c}{Validation} & \multicolumn{2}{c}{Held-out} & \multicolumn{2}{c}{Validation} & \multicolumn{2}{c}{Held-out} \\
    \cmidrule(lr){2-3}\cmidrule(lr){4-5}\cmidrule(lr){6-7}\cmidrule(l){8-9}
    & Stage~1 & Stage~2 & Stage~1 & Stage~2 & LR-aware & LR-blind & LR-aware & LR-blind \\
    \midrule
    Input tokens (thousands) & 1.9 & 88.4 & 2.2 & 86.4 & 136.2 & 121.6 & 163.8 & 146.4 \\
    Output tokens (thousands) & 1.5 & 20.8 & 1.6 & 17.8 & 278.4 & 262.7 & 332.4 & 327.2 \\
    Cost (USD) & \$0.03 & \$0.61 & \$0.04 & \$0.57 & \$5.95 & \$5.59 & \$7.11 & \$6.95 \\
    \bottomrule
  \end{tabular*}
\end{table}
\endgroup

Table~\ref{tab:short_llm_usage} averages per-run LLM usage over all LLM--LLM relaxation runs and all evaluated primal heuristic runs.
Costs use Claude Opus~5.5 API rates of \$4 and \$20 per million input and output tokens, plus \$5 and \$0.20 per million cache writes and reads, which are counted under input tokens.
For the Relaxation Design Agent, Stage~1 makes one call per initial pool and Stage~2 makes one formulation-summary call plus two calls per round, averaging 10.3 seconds per call.
For the Primal Heuristic Agent, each run makes one strategy call, three implementation calls, three revision calls, and one correction call per error. 
Since the six code calls return complete Julia programs, output tokens account for 68\% of all tokens, and the LR-aware and LR-blind variants differ in cost by under 7\%.

\section{Time-Budget Sensitivity of Dual AUC}
\label{sec:time_budget_sensitivity}

\begin{table}[H]
  \centering
  \caption{Time-budget sensitivity of relaxation selection comparing cutoffs $h$ against $T=120$ seconds under the same metric. Match is the percentage of instances within 0.001 of the full-budget best score, and Regret is the mean score loss relative to the best relaxation at $T$.}
  \label{tab:cutoff_stability}
  \footnotesize
  \setlength{\tabcolsep}{4pt}
  \begin{tabular}{@{}ll rr rr rr@{}}
    \toprule
    & & \multicolumn{2}{c}{\textbf{Validation}} & \multicolumn{2}{c}{\textbf{Held-out}} & \multicolumn{2}{c}{\textbf{Total}} \\
    \cmidrule(lr){3-4}\cmidrule(lr){5-6}\cmidrule(l){7-8}
    \textbf{Cutoff} & \textbf{Score} & \textbf{Match (\%)} & \textbf{Regret} & \textbf{Match (\%)} & \textbf{Regret} & \textbf{Match (\%)} & \textbf{Regret} \\
    \midrule
    \multirow{2}{*}{30\,s} & Dual AUC & \textbf{85.0} & \textbf{0.005} & \textbf{73.3} & \textbf{0.010} & \textbf{79.2} & \textbf{0.008} \\
     & Best Bound & 72.5 & 0.022 & 70.4 & 0.022 & 71.4 & 0.022 \\
    \midrule
    \multirow{2}{*}{60\,s} & Dual AUC & \textbf{95.0} & \textbf{0.001} & \textbf{96.7} & \textbf{0.001} & \textbf{95.8} & \textbf{0.001} \\
     & Best Bound & 87.2 & 0.015 & 77.0 & 0.009 & 82.1 & 0.012 \\
    \midrule
    \multirow{2}{*}{90\,s} & Dual AUC & \textbf{95.0} & \textbf{$<$0.001} & \textbf{98.3} & \textbf{$<$0.001} & \textbf{96.7} & \textbf{$<$0.001} \\
     & Best Bound & 90.4 & 0.011 & 83.5 & 0.008 & 87.0 & 0.010 \\
    \bottomrule
  \end{tabular}
\end{table}

Truncating all Phase~1 candidate trajectories at $h \in \{30,60,90\}$ seconds compares the time-budget sensitivity of Dual AUC against endpoint Best Bound $q_I(b_z(h))$.
Across both suites in \cref{tab:cutoff_stability}, Dual AUC achieves higher Match and lower Regret at every cutoff, exceeding 95\% Match with Regret at most 0.001 from 60 seconds onward.
Its smallest margin occurs at 30 seconds on the held-out suite, where slow-starting relaxations accumulate little early area.
Integrating over $[0,T]$ in \cref{eq:dual_auc} retains progress over $[0,h]$ with weight $h/T$ and dampens sensitivity to late bound jumps that can flip endpoint rankings, so Dual AUC rankings settle well before the 120-second budget.
\clearpage
\begin{landscape}
\section{Validation Set}
\label{app:short_validation}

\begingroup
  \footnotesize
  \setlength{\baselineskip}{9pt}
  \renewcommand{\arraystretch}{0.8}
  \setlength{\LTpre}{0pt}
  \setlength{\tabcolsep}{3pt}
  \setlength{\LTcapwidth}{\linewidth}
  \begin{longtable}{@{}>{\fontsize{9}{10}\selectfont\raggedright\arraybackslash}p{\dimexpr0.17\linewidth-2\tabcolsep\relax}>{\fontsize{9}{10}\selectfont\raggedright\arraybackslash}p{\dimexpr0.11\linewidth-2\tabcolsep\relax}>{\fontsize{9}{10}\selectfont\raggedleft\arraybackslash}p{\dimexpr0.075\linewidth-2\tabcolsep\relax}>{\fontsize{9}{10}\selectfont\raggedleft\arraybackslash}p{\dimexpr0.085\linewidth-2\tabcolsep\relax}>{\fontsize{9}{10}\selectfont\raggedright\arraybackslash}p{\dimexpr0.28\linewidth-2\tabcolsep\relax}>{\fontsize{9}{10}\selectfont\raggedright\arraybackslash}p{\dimexpr0.28\linewidth-2\tabcolsep\relax}@{}}
    \caption{Validation benchmarks. Original MILP gap is the range of Gurobi MIP gaps after 120 seconds on all cores over the three instances of each level; 0.0 denotes proven optimality.}
  \label{tab:short_validation}\\
    \toprule
    & & \multicolumn{2}{c}{\shortstack{\textbf{Original MILP}\\\textbf{gap (\%)}}} & \multicolumn{2}{c}{\textbf{Relaxation}} \\
    \cmidrule(lr){3-4}\cmidrule(l){5-6}
    \textbf{Problem} & \textbf{Source} & \textbf{Medium} & \textbf{Hard} & \textbf{Expert} & \textbf{AutoLR} \\
    \midrule
    \endfirsthead
    \multicolumn{6}{@{}l}{\textit{Validation benchmarks (continued)}}\\
    \toprule
    & & \multicolumn{2}{c}{\shortstack{\textbf{Original MILP}\\\textbf{gap (\%)}}} & \multicolumn{2}{c}{\textbf{Relaxation}} \\
    \cmidrule(lr){3-4}\cmidrule(l){5-6}
    \textbf{Problem} & \textbf{Source} & \textbf{Medium} & \textbf{Hard} & \textbf{Expert} & \textbf{AutoLR} \\
    \midrule
    \endhead
    \bottomrule
    \endlastfoot
    P01: Capacitated multi-item lot sizing &
    \citet{Gao2024Lagrange} &
    0.1--0.3 & 1.5--3.1 &
    Dualizes capacity: independent single-item lot-sizing subproblems. &
    Reproduces the expert relaxation. \\
    \addlinespace[3.5pt]
    P02: Constrained open-pit production scheduling &
    \citet{moosavi2014improvement} &
    0.0--3.8 & 11.0--26.0 &
    Dualizes grade and capacity limits: one block retaining precedence. &
    Reproduces the expert relaxation. \\
    \addlinespace[3.5pt]
    P03: Uncertainty-based open-pit production scheduling &
    \citet{tolouei2021optimisation} &
    1.8--2.5 & 4.7--8.2 &
    Dualizes grade and capacity limits: one block retaining precedence. &
    Dualizes grade limits and precedence: one precedence-free scheduling block retaining capacity limits. \\
    \addlinespace[3.5pt]
    P04: Green capacitated multi-item lot sizing with parallel machines &
    \citet{Wu2018GreenCapacitated} &
    1.3--2.0 & 4.4--4.8 &
    Dualizes capacity, emission limits, and machine-purchase constraints: a machine-purchase subproblem and single-item lot-sizing subproblems. &
    Dualizes capacity, emission limits, and one-machine-per-item constraints: joint lot-sizing and machine-purchase subproblems. \\
    \addlinespace[3.5pt]
    P05: Order batching in steelmaking and continuous casting &
    \citet{Tang2011OrderBatching} &
    1.8--2.3 & 3.7--4.0 &
    Dualizes cross-lot order assignment: one block, separating by lot only when charge counts are fixed. &
    Additionally dualizes minimum slab-weight bounds: same single-block structure. \\
    \addlinespace[3.5pt]
    P06: Capacity-constrained fast-fashion transshipment &
    \citet{NaderiKilicDasci2020} &
    1.9--2.1 & 8.2--10.9 &
    Dualizes sales-availability and route-activation constraints: closed-form, knapsack, and generalized-assignment subproblems. &
    Dualizes store outbound, destination, and route-activation constraints: independent per-product transshipment and sales subproblems. \\
    \addlinespace[3.5pt]
    P07: Two-stage hybrid flow shop with limited storage &
    \citet{GuMemarZinder2019} &
    0.0--14.6 & 16.1--20.4 &
    Dualizes machine and storage capacities: independent per-batch subproblems. &
    Additionally dualizes batch-start indicator constraints: same per-batch structure, with isolated batch-start indicators. \\
    \addlinespace[3.5pt]
    P08: Batch processing in a flexible flow shop &
    \citet{GuMemarZinder2018} &
    0.0 & 8.6--27.1 &
    Dualizes machine and storage capacities: independent per-batch subproblems. &
    Dualizes stage-2 machine capacity, storage capacity, and batch-scheduling constraints: one block with stage-1 machine capacity. \\
    \addlinespace[3.5pt]
    P09: Integrated production and maintenance planning &
    \citet{AlaouiSelsouli2012} &
    0.8--1.1 & 4.0--5.2 &
    Dualizes shared resource capacity: single-item lot-sizing subproblems with shortage costs and a maintenance-planning subproblem. &
    Additionally dualizes maintenance-timing constraints: same lot-sizing subproblems and a simpler maintenance-planning subproblem. \\
    \addlinespace[3.5pt]
    P10: Hybrid flow shop scheduling with earliness and tardiness &
    \citet{DingLiuZhuLi2014} &
    0.0--5.0 & 25.9--29.6 &
    Dualizes stage precedence: stage-occupancy constraints still couple stages into one block. &
    Dualizes machine capacity: independent per-job multi-stage scheduling subproblems. \\
  \end{longtable}
\endgroup
\end{landscape}

\begin{landscape}
\section{Held-Out Set}
\label{app:short_heldout}

\begingroup
  \footnotesize
  \setlength{\baselineskip}{9pt}
  \renewcommand{\arraystretch}{0.8}
  \setlength{\LTpre}{0pt}
  \setlength{\LTpost}{4pt}
  \setlength{\tabcolsep}{3pt}
  \setlength{\LTcapwidth}{\linewidth}
  \begin{longtable}{@{}>{\fontsize{9}{10}\selectfont\raggedright\arraybackslash}p{\dimexpr0.17\linewidth-2\tabcolsep\relax}>{\fontsize{9}{10}\selectfont\raggedright\arraybackslash}p{\dimexpr0.11\linewidth-2\tabcolsep\relax}>{\fontsize{9}{10}\selectfont\raggedleft\arraybackslash}p{\dimexpr0.075\linewidth-2\tabcolsep\relax}>{\fontsize{9}{10}\selectfont\raggedleft\arraybackslash}p{\dimexpr0.085\linewidth-2\tabcolsep\relax}>{\fontsize{9}{10}\selectfont\raggedright\arraybackslash}p{\dimexpr0.28\linewidth-2\tabcolsep\relax}>{\fontsize{9}{10}\selectfont\raggedright\arraybackslash}p{\dimexpr0.28\linewidth-2\tabcolsep\relax}@{}}
    \caption{Held-out benchmarks. Original MILP gap follows \cref{tab:short_validation}. $^{\dagger}$No majority among the nine medium runs; selected by highest mean Dual AUC as in \cref{sec:autolr_selection}.}
  \label{tab:short_heldout}\\
    \toprule
    & & \multicolumn{2}{c}{\shortstack{\textbf{Original MILP}\\\textbf{gap (\%)}}} & \multicolumn{2}{c}{\textbf{Relaxation}} \\
    \cmidrule(lr){3-4}\cmidrule(l){5-6}
    \textbf{Problem} & \textbf{Source} & \textbf{Medium} & \textbf{Hard} & \textbf{hMETIS} & \textbf{AutoLR} \\
    \midrule
    \endfirsthead
    \multicolumn{6}{@{}l}{\textit{Held-out benchmarks (continued)}}\\
    \toprule
    & & \multicolumn{2}{c}{\shortstack{\textbf{Original MILP}\\\textbf{gap (\%)}}} & \multicolumn{2}{c}{\textbf{Relaxation}} \\
    \cmidrule(lr){3-4}\cmidrule(l){5-6}
    \textbf{Problem} & \textbf{Source} & \textbf{Medium} & \textbf{Hard} & \textbf{hMETIS} & \textbf{AutoLR} \\
    \midrule
    \endhead
    \midrule
    \multicolumn{6}{r@{}}{\textit{Continued on the next page}}\\
    \endfoot
    \bottomrule
    \endlastfoot
    T01: Transportation planning with rules &
    MIPLIB-NL &
    0.0--2.7 & 4.2--6.4 &
    Dualizes share limits plus coverage and policy subsets: multiple blocks. &
    Dualizes mutual-exclusion constraints: one block. \\
    \addlinespace[2.5pt]
    T02$^{\dagger}$: Diagram entity placement with anchoring &
    MIPLIB-NL &
    0.0 & 3.8--11.4 &
    Dualizes anchor-deviation, relative-position, and Manhattan bounds: one row-axis and one column-axis block. &
    Dualizes column gaps, column-order consistency, and Manhattan bounds: one block coupling both axes. \\
    \addlinespace[2.5pt]
    T03: Robust linear regression with outlier exclusion &
    MIPLIB-NL &
    0.0 & 4.4--7.7 &
    Dualizes outlier budgets plus subsets of fitting and big-M rows: the subproblem becomes unbounded. &
    No relaxation selected: big-M dualization is unbounded, while retention yields negative bounds via $M$. \\
    \addlinespace[2.5pt]
    T04: Molecular redesign as state and interaction selection &
    MIPLIB-NL &
    0.0 & 5.6--7.0 &
    Dualizes pair-selection and both endpoint constraints: one block per site, with pairwise interactions unconstrained. &
    Dualizes left-endpoint constraints: one block per site, with pairwise interactions attached to their right endpoint. \\
    \addlinespace[2.5pt]
    T05: Residential tower design with sector and ownership quotas &
    \citet{ferchtandiker2025efficient} &
    0.0--1.0 & 2.5--4.7 &
    Dualizes tower-wide quotas and the layout constraint for one floor: per-floor blocks, with that floor split into small blocks. &
    Dualizes minimum sector-share quotas and owner--sector eligibility constraints: one tower-wide block. \\
    \addlinespace[2.5pt]
    T06$^{\dagger}$: Two-echelon truck and drone routing for post-disaster relief &
    \citet{nguyen2026twoechelon} &
    0.0 & 41.8--43.0 &
    Dualizes truck capacity, cluster coverage, and timing subsets: separate assignment and routing blocks. &
    Dualizes truck capacity and truck--drone load linearization constraints: one block. \\
    \addlinespace[2.5pt]
    T07: Integrated patient-room and nurse assignment &
    \citet{brandt2025integrated} &
    0.0--3.7 & 13.8--14.0 &
    Dualizes subsets of cross-day continuity and fairness rows: multi-day blocks retaining room--nurse coupling. &
    Dualizes room--nurse assignment constraints: one full-horizon nurse block, one room block, and per-shift walking blocks. \\
    \addlinespace[2.5pt]
    T08: Language-model placement on heterogeneous GPUs &
    \citet{mei2025helix} &
    0.0 & 288.9--309.5 &
    Dualizes flow conservation plus throughput and layer-covering subsets: two balanced blocks. &
    Dualizes both layer-covering conditions: one block. \\
    \addlinespace[2.5pt]
    T09: Adaptive rescheduling of patients in operating-room blocks &
    \citet{kamran2020column} &
    0.0--5.7 & 24.9--83.3 &
    Dualizes subsets of assignment and sequencing rows: two blocks retaining block sequencing. &
    Dualizes within-block big-M sequencing constraints: one unordered assignment block. \\
    \addlinespace[2.5pt]
    T10$^{\dagger}$: Cross-docking network scheduling and loading &
    \citet{castellucci2021network} &
    0.0--3.2 & 5.1--10.7 &
    Dualizes subsets of routing, timing, destination, and box-separation rows: multiple blocks retaining packing. &
    Dualizes single-destination and box-separation constraints: one routing block and per-axis packing blocks. \\
  \end{longtable}
  \noindent\parbox{\linewidth}{\fontsize{9}{10}\selectfont
  MIPLIB-NL~\citep{li2026constructing}: T01 = \texttt{bab1}, T02 = \texttt{graphdraw-mainerd}, T03 = \texttt{neos-3754480-nidda}, T04 = \texttt{proteindesign122trx11p8}.\\
  LLMs-For-Optimization-Reformulations~\citep{ferchtandiker2025efficient}: T05 = \texttt{dutch\_housing\_problem}.}
\endgroup
\end{landscape}

\end{document}